\documentclass{article}

\PassOptionsToPackage{numbers,sort&compress}{natbib}
\usepackage[final,main]{neurips_2026}

\usepackage[utf8]{inputenc}
\usepackage[T1]{fontenc}

\usepackage{url}
\usepackage{booktabs}
\usepackage{array}
\usepackage{multirow}
\usepackage{amsfonts}
\usepackage{amsmath}
\usepackage{amssymb}
\usepackage{amsthm}
\usepackage{tikz}
\usetikzlibrary{arrows.meta, positioning}
\usepackage{mathtools}
\usepackage{nicefrac}
\usepackage{microtype}
\usepackage{xcolor}
\usepackage{bm}
\usepackage{enumitem}
\usepackage{algorithm}
\usepackage{algorithmic}
\usepackage{placeins}
\usepackage{graphicx}
\usepackage{hyperref}

\theoremstyle{plain}
\newtheorem{theorem}{Theorem}[section]
\newtheorem{proposition}[theorem]{Proposition}
\newtheorem{lemma}[theorem]{Lemma}

\theoremstyle{definition}
\newtheorem{definition}[theorem]{Definition}
\newtheorem{assumption}[theorem]{Assumption}

\newtheorem{remark}[theorem]{Remark}

\newcommand{\E}{\mathbb{E}}
\newcommand{\R}{\mathbb{R}}
\newcommand{\N}{\mathbb{N}}
\newcommand{\Prob}{\mathbb{P}}
\newcommand{\Var}{\mathrm{Var}}

\newcommand{\cF}{\mathcal{F}}
\newcommand{\cT}{\mathcal{T}}

\newcommand{\VaR}{\mathrm{VaR}}
\newcommand{\CVaR}{\mathrm{CVaR}}

\newcommand{\ind}{\mathbf{1}}

\title{CVaR Deep Optimal Stopping}

\author{%
  Yechan Jeong \\
  Seoul National University \\
  \texttt{yechan@snu.ac.kr}
  \And
  Youngjoon Hong \\
  Seoul National University \\
  \texttt{hongyj@snu.ac.kr}
}

\begin{document}
\raggedbottom

\maketitle

\begin{abstract}
Risk-sensitive optimal stopping is useful when poor stopping outcomes carry substantial costs, as in option exercise, sequential medical treatment, and safety-critical control. A natural objective is Conditional Value at Risk (CVaR) of the stopped reward. We develop a theoretical framework and algorithm for CVaR deep optimal stopping. The framework and algorithm address four difficulties in applying risk-sensitive reinforcement learning to optimal stopping: state augmentation, score-function pathology, blindness to success, and dependence on auxiliary networks. We isolate three structural properties of optimal stopping: single-point reward, action-independent dynamics, and sum-zero gradient. We exploit these properties to establish Markov optimality without state augmentation and derive a pathwise gradient estimator. Under bounded policy derivatives, the pathwise estimator has bounded second moments near deterministic policies, avoiding the divergence that can occur with score-function estimators. Building on this framework, we propose DIOS (Distribution-Informed Optimal Stopping) to address the remaining two difficulties. A per-time-step soft penalty combined with tail-level annealing mitigates blindness to success. Estimating the reward quantile directly from each batch removes the need for an auxiliary critic. For a fixed tail level, we establish a finite-time stationarity guarantee under regularity assumptions along the training iterates. Across the 24 main Bermudan max-call configurations, DIOS achieves the highest mean CVaR and the shortest recorded training time among the tested CVaR baselines.
\end{abstract}

\section{Introduction}
\label{sec:intro}

\paragraph{Risk-sensitive deep optimal stopping.}
Risk-sensitive optimal stopping is useful when poor stopping outcomes carry substantial costs, as in option exercise, sequential medical treatment, and safety-critical control. Maximizing expected reward alone may give insufficient weight to adverse stopping outcomes. We study the Conditional Value at Risk (CVaR) of the stopped reward, which measures the average reward over a chosen fraction of the worst outcomes~\citep{rockafellar2000optimization,rockafellar2002conditional}.

\paragraph{Coherent risk measure and CVaR.}
A coherent risk measure satisfies monotonicity, translation invariance, positive homogeneity, and subadditivity~\citep{artzner1999}. The negative of lower-tail CVaR is a coherent risk measure for rewards~\citep{rockafellar2002conditional}. We choose CVaR to optimize the average reward over a prescribed fraction of the worst stopping outcomes. The Rockafellar--Uryasev representation decomposes CVaR optimal stopping into an outer optimization over a scalar variable and an inner risk-neutral stopping problem. Standard stopping theory and deep optimal stopping methods can therefore be applied to the inner stopping problem.

\paragraph{Deep optimal stopping.}
The Snell envelope characterizes risk-neutral optimal stopping~\citep{shiryaev2008optimal}, and the regression-based approach of~\citet{longstaff2001valuing} provides a method for American option pricing. Deep learning is used to address the curse of dimensionality in optimal stopping. \citet{becker2019deep} introduced \emph{deep optimal stopping}, a backward-induction method that parameterizes each per-step stopping probability with a neural network and trains the networks in backward order. Neural-network regression methods approximate continuation values~\citep{lapeyre2021neural,herrera2024optimal}. Other methods jointly optimize stopping probabilities across time~\citep{becker2021solving,venkata2023}, with shared policy parameters in~\citet{venkata2023}. Language models have also been applied to stopping decisions in sales conversations~\citep{manzoor2025}. We study CVaR deep optimal stopping, replacing expected-reward maximization with lower-tail CVaR maximization.

\paragraph{Risk-sensitive reinforcement learning.}
For CVaR objectives and constraints, \citet{bauerle2011}, \citet{chow2014}, and \citet{chow2018risk} obtain Bellman recursions by augmenting the state with risk or cumulative-reward variables. \citet{tamar2015cvar} derive a score-function CVaR policy-gradient (CVaR-PG) estimator. Methods for improving sample efficiency in CVaR policy optimization modify sampling, policy parameterization, or gradient estimation. \citet{greenberg2022} pair a soft-risk schedule with a cross-entropy environment sampler. \citet{luo2024mixture} combine a risk-averse policy with a risk-neutral policy trained via Implicit Q-Learning~\citep{kostrikov2022}. \citet{kim2024pcvar} weight individual costs by predicted tail probabilities from a separate network. \citet{mead2025retcap} compute policy gradients from trajectory returns capped at a value updated during training. \citet{luo2026cvarpg} mix the CVaR-PG gradient with a quantile-gradient term computed by an auxiliary quantile critic. These CVaR methods can be applied to optimal stopping, a finite-horizon MDP with a binary stop/continue action and action-independent state process. Applying risk-sensitive reinforcement learning (RL) to optimal stopping raises four difficulties.

\paragraph{Four structural problems of risk-sensitive RL.}
\begin{itemize}[leftmargin=2em,itemsep=2pt,topsep=2pt]
\item[\textbf{(P1)}] \textbf{State augmentation.} Static CVaR objectives can be time-inconsistent~\citep{shapiro2009time}. Existing CVaR-MDP methods obtain Bellman recursions by adding variables for accumulated costs to the state~\citep{bauerle2011,chow2014,chow2018risk}.
\item[\textbf{(P2)}] \textbf{Score-function pathology.} Near deterministic policies, the second moment of the trajectory-level score-function estimator can diverge or vanish, depending on the policy parameterization. Several methods retain a score-function CVaR-PG component~\citep{greenberg2022,luo2024mixture,luo2026cvarpg}.
\item[\textbf{(P3)}] \textbf{Blindness to success.} The CVaR-PG gradient estimator assigns zero weight to trajectories with rewards above the estimated quantile. The CVaR-PG gradient estimate can also vanish when the sampled lower-tail rewards coincide~\citep{greenberg2022,luo2026cvarpg}.
\item[\textbf{(P4)}] \textbf{Dependence on auxiliary networks.} Several methods use learned tail probabilities, values, or quantiles in policy updates~\citep{kim2024pcvar,luo2024mixture,luo2026cvarpg}. The auxiliary networks are trained on data generated by the policy.
\end{itemize}

Optimal stopping addresses (P1) and permits a pathwise treatment of (P2). DIOS addresses (P3) and removes the auxiliary critic in (P4).

\paragraph{Framework: addressing (P1) and (P2).}
The single-point reward (S1) and the envelope reduction give Markov optimality without state augmentation. Action-independent dynamics (S2) permit differentiation of the conditional stopping probabilities without a trajectory score term. The sum-zero gradient (S3) makes the pathwise gradient invariant to a common shift of the transformed rewards. Under the policy-derivative bounds, pathwise gradient moments remain bounded near deterministic policies.

\paragraph{DIOS and convergence: addressing (P3) and (P4).}
We propose \emph{Distribution-Informed Optimal Stopping} (DIOS), a pathwise algorithm for CVaR optimal stopping. A per-time-step soft penalty allows rewards at other stopping times to contribute even when the sampled stopping reward exceeds the estimated quantile. Training starts with a tail level of $1$, then gradually reduces the tail level to its target. Direct quantile estimation from each batch removes the auxiliary critic. Theorem~\ref{thm:conv} gives a finite-time best-iterate stationarity guarantee when the analyzed iterates remain in a domain satisfying the regularity assumptions.

\paragraph{Contributions.}
\begin{itemize}[leftmargin=2em,itemsep=2pt,topsep=2pt,beginpenalty=10000]
\item We exploit (S1)--(S3) to develop a CVaR optimal stopping framework addressing (P1) and (P2). We establish the envelope reduction (Proposition~\ref{thm:envelope}) and prove Markov optimality without state augmentation (Theorem~\ref{thm:markov}).
\item We derive a pathwise gradient estimator (Lemma~\ref{thm:pathwise}) and bound its second moments under the policy-derivative and reward-moment assumptions (Lemma~\ref{thm:variance}). We compare pathwise and score-function gradients near deterministic policies (Theorem~\ref{cor:pathwise-vs-score}).
\item We propose \emph{Distribution-Informed Optimal Stopping} (DIOS, Algorithm~\ref{alg:dios}). DIOS combines a per-time-step soft penalty with tail-level annealing to mitigate (P3). Direct quantile estimation removes the auxiliary critic in (P4). At a fixed tail level, we establish a finite-time best-iterate stationarity guarantee under the regularity and parameter-domain conditions of Theorem~\ref{thm:conv}.
\item We demonstrate across 24 Bermudan max-call configurations that DIOS achieves the highest mean CVaR and the shortest recorded training time among the tested CVaR baselines~\citep{tamar2015cvar,kim2024pcvar,luo2024mixture,mead2025retcap,luo2026cvarpg}.
\end{itemize}

\section{Deep optimal stopping and CVaR}
\label{sec:setup}
\label{sec:dos}
\subsection{Deep optimal stopping}
\label{sec:dos-policy}

\paragraph{Setup.}
We work in the finite-horizon, discrete-time, model-free setting. Fix a probability space $(\Omega, \cF, P)$ and a finite horizon $T \in \N$. Let $(S_t)_{t=0}^T$ be an $\R^d$-valued Markov process and $(\cF_t)_{t=0}^T$ its natural filtration, $\cF_t := \sigma(S_0, \dots, S_t)$. Define $R_t := g_t(S_t)$, where $g_t : \R^d \to \R_+$ is Borel-measurable. The joint distribution of $(S_0,\ldots,S_T)$ is unknown and accessed only through samples. Appendix~\ref{app:exp-verif} discusses the assumptions for the Bermudan max-call experiments in Section~\ref{sec:exp}.

\begin{assumption}[Integrability]
\label{asm:reg}
$R_t \in L^1(P)$ for every $t \in \{0, \dots, T\}$.
\end{assumption}

\paragraph{The optimal stopping problem.}
Let $\cT$ be the set of random variables $\tau:\Omega\to\{0,\ldots,T\}$ satisfying $\{\tau=t\}\in\cF_t$ for every $t$. 
The optimal stopping problem is to find $\tau^* \in \cT$ attaining
\begin{equation}
\label{eq:os}
V^* \;:=\; \sup_{\tau \in \cT} \E\!\left[R_\tau\right].
\end{equation}
Existence of $\tau^*$ under Assumption~\ref{asm:reg} is classical~\citep{shiryaev2008optimal}. Define the Snell envelope $(V_t)_{t=0}^T$ by backward induction:
\begin{equation}
\label{eq:snell}
V_T := R_T, \quad V_t := \max\!\left(R_t,\ \E[V_{t+1} \mid \cF_t]\right) \quad \text{for } t < T.
\end{equation}
Assumption~\ref{asm:reg} ensures that $V_t$ is integrable for every $t=0,\ldots,T$. For $t<T$, an optimal rule stops when the immediate reward $R_t$ satisfies $R_t\ge\E[V_{t+1}\mid\cF_t]$ and continues otherwise, with a forced stop at $T$. The supremum in~\eqref{eq:os} equals $\E[V_0]$, attained at $\tau^* := \inf\{t \ge 0 : V_t = R_t\}$.

\paragraph{Stochastic stopping policies~\citep{venkata2023}.}
Let $\phi=(\phi_0,\ldots,\phi_T)$ consist of Borel functions $\phi_t:\R^d\to[0,1]$ with $\phi_T\equiv1$. At time $t$, conditional on not having stopped earlier, the policy stops with probability $\phi_t(S_t)$. Let $U_0,\ldots,U_T$ be independent uniform random variables on $[0,1]$, independent of the state process. On the enlarged probability space, define
\begin{equation}
\label{eq:tau}
\tau \;:=\; \inf\!\left\{t \,:\, U_t \le \phi_t(S_t)\right\}.
\end{equation}
Conditioning on $\cF_T$ fixes $S_0,\ldots,S_T$, leaving only the randomness in $U_0,\ldots,U_T$. Define $\psi_t:=P(\tau=t\mid\cF_T)$ and $\Psi_t:=\prod_{s<t}(1-\phi_s(S_s))$, with $\Psi_0:=1$. Then
\begin{equation}
\label{eq:psi}
\psi_t \;=\; \phi_t(S_t) \cdot \prod_{s < t}\!\left(1 - \phi_s(S_s)\right) \;=\; \Psi_t \, \phi_t(S_t), \qquad t = 0, \dots, T.
\end{equation}
The factor $\Psi_t=P(\tau\ge t\mid\cF_T)$ gives the conditional probability of not stopping before time $t$. Since $\phi_T\equiv1$, the probabilities $\psi_t$ satisfy
\begin{equation}
\label{eq:psi-sum}
\sum_{t=0}^T \psi_t \;=\; 1 \qquad \text{pointwise on } \Omega.
\end{equation}
By the tower property,
\begin{equation}
\label{eq:pathwise}
\E\!\left[R_\tau\right] \;=\; \E\!\left[\E[R_\tau \mid \cF_T]\right] \;=\; \E\!\left[\sum_{t=0}^T P(\tau = t \mid \cF_T)\, R_t\right] \;=\; \E\!\left[\sum_{t=0}^T \psi_t\, R_t\right].
\end{equation}

\paragraph{Three structural properties of optimal stopping.}
We exploit the following three properties of the stopping formulation.

\begin{itemize}[leftmargin=2em,itemsep=2pt,topsep=2pt]
\item[\textbf{(S1)}] \emph{Single-point reward.} The reward $R_\tau = g_\tau(S_\tau)$ is realized at the single time index $\tau$, in contrast to cumulative-return MDPs where rewards $\sum_t r_t$ are summed over the trajectory.

\item[\textbf{(S2)}] \emph{Action-independent dynamics.} The state process $(S_t)_{t=0}^T$ evolves under a transition kernel that does not depend on the stopping action. 

\item[\textbf{(S3)}] \emph{Sum-zero gradient.} Parameterize the policy by $\theta\in\Theta$, where $\Theta\subset\R^p$ is open, and assume each $\theta \mapsto \psi_t^\theta$ is differentiable $P$-a.s. Since~\eqref{eq:psi-sum} holds for every $\theta$, differentiating in $\theta$ yields
\[
\sum_{t=0}^T \nabla_\theta \psi_t^\theta \;=\; 0 \qquad P\text{-a.s.}
\]
\end{itemize}

\subsection{CVaR optimal stopping}
Fix $\beta\in(0,1)$. For $X\in L^1(P)$, define the lower $\beta$-quantile by $\VaR_\beta(X) := \inf\{x : P(X \le x) \ge \beta\}$. Define the \emph{Conditional Value at Risk} by
\begin{equation}
\label{eq:cvar}
\CVaR_\beta(X) \;:=\; \tfrac{1}{\beta}\int_0^\beta \VaR_u(X)\, du.
\end{equation}
Define $\CVaR_1(X):=\E[X]$. Smaller $\beta$ focuses the objective on a smaller fraction of the lowest outcomes. The Rockafellar-Uryasev variational form expresses CVaR as a one-dimensional supremum~\citep{rockafellar2002conditional},
\begin{equation}
\label{eq:ru}
\CVaR_\beta(X) \;=\; \sup_{\nu \in \R} \Bigl\{\nu \;-\; \tfrac{1}{\beta}\, \E[(\nu - X)^+]\Bigr\},
\end{equation}
attained at $\nu^*=\VaR_\beta(X)$ for $0<\beta<1$. The quantity $(\nu-X)^+=\max\{\nu-X,0\}$ measures the shortfall below $\nu$. Replacing the expectation in~\eqref{eq:os} with CVaR yields the \emph{CVaR optimal stopping value}
\begin{equation}
\label{eq:cvar-os}
V^*_\beta \;:=\; \sup_{\tau \in \cT}\, \CVaR_\beta(R_\tau).
\end{equation}

\section{CVaR within the deep optimal stopping framework}
\label{sec:dios}

Building on the Rockafellar--Uryasev envelope reduction, we establish Markov optimality without state augmentation and derive a pathwise gradient formula.

\subsection{Envelope reduction}
\label{sec:envelope}

Proofs are given in Appendix~\ref{app:proofs}.

\begin{proposition}[CVaR envelope reduction]
\label{thm:envelope}
Under Assumption~\ref{asm:reg} and $\beta \in (0, 1)$,
\begin{equation}
\label{eq:envelope}
V^*_\beta \;=\; \sup_{\nu \in \R} \Bigl\{\nu \;+\; \tfrac{1}{\beta}\, \sup_{\tau \in \cT} \E[-(\nu - R_\tau)^+]\Bigr\},
\end{equation}
and there exist $\nu^* \in \R$ and $\tau^* \in \cT$ with $\nu^* = \VaR_\beta(R_{\tau^*})$ jointly realizing the supremum.
\end{proposition}

Proposition~\ref{thm:envelope} separates an outer maximization over the scalar $\nu$ from an inner risk-neutral stopping problem at fixed $\nu$. Define the per-step \emph{soft penalty} $\tilde R^\nu_t := -\tfrac{1}{\beta}(\nu - R_t)^+$. For the parameterized stopping policies $\phi^\theta$, define the objective for the inner maximization over $\theta$ by
\begin{equation}
\label{eq:Jnu}
J^{(\nu)}(\theta) \;:=\; \E\!\left[\sum_{t=0}^T \psi_t^\theta \cdot \tilde R^\nu_t\right].
\end{equation}
For fixed $\nu$, the Snell-envelope recursion and stopping-policy representation in Section~\ref{sec:dos-policy} apply with $R_t$ replaced by $\tilde R_t^\nu$. Appendix~\ref{app:oce} extends the outer--inner decomposition to an optimized certainty equivalent objective.

\subsection{Markov optimality without augmentation}
\label{sec:markov}

Applying CVaR recursively at each time generally changes the static CVaR objective~\citep{shapiro2009time}. CVaR-MDP methods instead obtain Bellman recursions by adding variables for accumulated costs to the state~\citep{bauerle2011,chow2014,chow2018risk}.

\begin{definition}[Markov stopping time]
\label{def:markov-stopping}
Define a stopping time $\tau\in\cT$ to be \emph{Markov} if there exist Borel sets $D_0,\ldots,D_T\subseteq\R^d$ with $D_T=\R^d$ such that $\tau=\inf\{t\ge0:S_t\in D_t\}$ $P$-a.s.
\end{definition}

\begin{theorem}[Markov optimality without augmentation]
\label{thm:markov}
Under Assumption~\ref{asm:reg} and $\beta\in(0,1)$, there exist $\nu^*\in\R$ and a Markov stopping time $\tau^*\in\cT$ such that
\[
V^*_\beta \;=\; \nu^* \;+\; \tfrac{1}{\beta}\, \E\!\left[-(\nu^* - R_{\tau^*})^+\right].
\]
\end{theorem}

\begin{remark}
For fixed $\nu^*$, (S1) makes the transformed reward $-(\nu^*-g_t(S_t))^+$ a function of the current state and time. By the Markov property, each conditional expectation in the Snell-envelope recursion is a function of $t$ and $S_t$, with no additional dependence on $S_0,\ldots,S_{t-1}$.
\end{remark}

\subsection{Pathwise and score-function gradients}\label{sec:pathwise}
For fixed $\nu$, we show that pathwise gradient second moments remain uniformly bounded under bounded policy derivatives, whereas score-function second moments can diverge as stopping probabilities approach zero.

By (S2), the joint distribution of $(S_0,\ldots,S_T)$ is independent of $\theta$. At fixed $\nu$, only $\psi_t^\theta$ in the integrand of~\eqref{eq:Jnu} depends on $\theta$.

\begin{assumption}[Differentiable policy with bounded gradients]
\label{asm:smooth}
Let $\|\cdot\|$ denote the Euclidean norm on $\R^p$. On an event of probability one, each map $\theta\mapsto\phi_t^\theta(S_t)$, $0\le t\le T$, is differentiable on the open domain $\Theta$. There exist $G,C_\psi\ge0$ such that, $P$-a.s.,
\[
\sup_{\theta\in\Theta}\max_{0\le t\le T}\|\nabla_\theta\phi_t^\theta(S_t)\|\le G,
\qquad
\sup_{\theta\in\Theta}\max_{0\le t\le T}\|\nabla_\theta\psi_t^\theta\|\le C_\psi.
\]
For every $\theta\in\Theta$ and $t<T$, $\phi_t^\theta(S_t)\in(0,1)$ $P$-a.s. The bound $C_\psi$ may depend on $T$.
\end{assumption}
Assumptions~\ref{asm:reg} and~\ref{asm:smooth} justify differentiation under the expectation in~\eqref{eq:Jnu}. Equation~\eqref{eq:Cpsi-loose} in Appendix~\ref{app:exp-verif} shows that $C_\psi=G(T+1)$ is a valid choice.

Define $\mathrm{sg}[x]$ to have value $x$ and zero derivative in automatic differentiation. Following~\citet{venkata2023}, define the surrogate loss for the transformed rewards $\tilde R_t^\nu$ at fixed $\nu\in\R$ by
\begin{equation}
\label{eq:dios-loss}
L(\theta) \;:=\; -\E\!\left[\sum_{t=0}^T \mathrm{sg}\!\left[\tilde R^\nu_t\right] \cdot \mathrm{sg}\!\left[\psi_t^\theta\right] \cdot \log \psi_t^\theta\right].
\end{equation}
Let $\nabla_\theta L$ denote the expected autograd gradient with $S_0,\ldots,S_T$ and $\nu$ held fixed. The stop-gradient factors restrict differentiation to $\log\psi_t^\theta$, giving $\psi_t^\theta\nabla_\theta\log\psi_t^\theta=\nabla_\theta\psi_t^\theta$.

\begin{lemma}[Pathwise gradient formula]\label{thm:pathwise}
Under Assumptions~\ref{asm:reg} and~\ref{asm:smooth}, for every fixed $\nu\in\R$,
\begin{equation}\label{eq:pathwise-grad}
\nabla_\theta J^{(\nu)}(\theta)
=\E\!\left[\sum_{t=0}^{T}\tilde R_t^\nu\nabla_\theta\psi_t^\theta\right]
=-\nabla_\theta L(\theta).
\end{equation}
\end{lemma}

\begin{remark}\label{rem:boundary}
Differentiating $\psi_t^\theta=\phi_t^\theta(S_t)\prod_{s<t}(1-\phi_s^\theta(S_s))$ gives
\[
\nabla_\theta\psi_t^\theta
=\Psi_t^\theta\nabla_\theta\phi_t^\theta(S_t)
-\phi_t^\theta(S_t)\sum_{s<t}
\left(\prod_{\substack{u<t\\u\ne s}}(1-\phi_u^\theta(S_u))\right)
\nabla_\theta\phi_s^\theta(S_s).
\]
No term divides by $\phi_s^\theta(S_s)$ or $1-\phi_s^\theta(S_s)$. Thus $\nabla_\theta\psi_t^\theta$ remains defined at probabilities $0$ or $1$ whenever $\nabla_\theta\phi_s^\theta(S_s)$ exists for every $s\le t$.
\end{remark}
\begin{assumption}[$L^2$ payoff]
\label{asm:l2}
Assume $R_t\in L^2(P)$ for $0\le t\le T$, and define $M:=\max_{0\le t\le T}\|R_t\|_{L^2(P)}$.
\end{assumption}

Let $(S_t^{(i)})_{t=0}^T$, $i=1,\ldots,B$, be independent copies of $(S_t)_{t=0}^T$. Let the superscript $(i)$ denote evaluation on trajectory $i$. Define the empirical gradient by
\[
\hat\nabla L(\theta):=-\frac1B\sum_{i=1}^B\sum_{t=0}^T
\tilde R_t^{\nu,(i)}\nabla_\theta\psi_t^{\theta,(i)}.
\]
By (S3), subtracting the mean $\frac{1}{T+1}\sum_{s=0}^T\tilde R_s^{\nu,(i)}$ from each $\tilde R_t^{\nu,(i)}$ leaves $\hat\nabla L(\theta)$ unchanged.
Let $\tau_\theta$ be given by~\eqref{eq:tau} for the policy $\phi^\theta$, and define $\nu^*(\theta):=\VaR_\beta(R_{\tau_\theta})$.
\begin{lemma}[Second moment and variance]\label{thm:variance}
Under Assumptions~\ref{asm:reg},~\ref{asm:smooth}, and~\ref{asm:l2}, fix $\theta\in\Theta$ and $\nu\in\R$. Define $\tilde M:=\max_t\|\tilde R_t^\nu\|_{L^2(P)}$.
\begin{enumerate}[label=\textbf{(\alph*)},leftmargin=2em,itemsep=1pt,topsep=1pt]
\item For every $B\ge1$,
\[
\E[\|\hat\nabla L(\theta)\|^2]\le C_\psi^2(T+1)^2\tilde M^2.
\]
\item For $B$ independent state trajectories,
\[
\E[\|\hat\nabla L(\theta)-\nabla L(\theta)\|^2]
\le \frac{C_\psi^2(T+1)^2\tilde M^2}{B}.
\]
\item If $M=0$, define $c_\beta:=0$. If $M>0$, define
\[
c_\beta:=\sup_{\vartheta\in\Theta}\max_{0\le t\le T}\frac{\E[((\nu^*(\vartheta)-R_t)^+)^2]}{M^2}.
\]
Then $c_\beta\le(T+1)/(1-\beta)$. At $\nu=\nu^*(\theta)$,
\begin{equation}\label{eq:variance-bound}
\E[\|\hat\nabla L(\theta)-\nabla L(\theta)\|^2]
\le\frac{4C_\psi^2c_\beta T^2M^2}{\beta^2B}.
\end{equation}
\end{enumerate}
\end{lemma}
In~\eqref{eq:variance-bound}, the factor $C_\psi^2T^2$ becomes $G^2T^2(T+1)^2$ for the choice $C_\psi=G(T+1)$. If $C_\psi$ is independent of $T$, $C_\psi^2T^2$ is proportional to $T^2$.

Let $\tau=\tau_\theta$. For $t<T$, let $a_t=\ind\{U_t\le\phi_t^\theta(S_t)\}$ and define
$\pi_\theta(a_t\mid S_t):=(\phi_t^\theta(S_t))^{a_t}(1-\phi_t^\theta(S_t))^{1-a_t}$.
Define the likelihood of the sampled actions up to $\tau$, conditional on $(S_0,\ldots,S_T)$, by
$\pi_\theta(\tau):=\prod_{t=0}^{\min\{\tau,T-1\}}\pi_\theta(a_t\mid S_t)$.
For comparison with existing score-function CVaR policy-gradient methods~\citep{tamar2015cvar}, define the loss for automatic differentiation at fixed $\nu$ by
\[
L_{\rm score}(\theta) \;:=\; -\E\!\left[\mathrm{sg}\!\left[\tilde R^\nu_\tau\right] \cdot \log \pi_\theta(\tau)\right].
\]
With sampled states, actions, and $\tau$ held fixed, define
\[
\hat\nabla L_{\rm score}(\theta):=-\tilde R^\nu_\tau\nabla_\theta\log\pi_\theta(\tau),
\qquad
\nabla_\theta L_{\rm score}(\theta):=\E[\hat\nabla L_{\rm score}(\theta)].
\]
Under Assumptions~\ref{asm:reg} and~\ref{asm:smooth}, $-\nabla_\theta L_{\rm score}(\theta)=\nabla_\theta J^{(\nu)}(\theta)$.
For the single-trajectory comparison $B=1$, $\hat\nabla L(\theta)$ sums over all candidate stopping times, whereas $\hat\nabla L_{\rm score}(\theta)$ uses the reward at the sampled stopping time $\tau$.
\begin{theorem}[Bounded pathwise second moments and score-function divergence]
\label{cor:pathwise-vs-score}
Under Assumptions~\ref{asm:reg},~\ref{asm:smooth}, and~\ref{asm:l2}, fix $\nu\in\R$ and a sequence $(\theta_n)\subset\Theta$. Let $B=1$. Suppose there exist $t^*\in\{0,\dots,T-1\}$, $\epsilon_0>0$, and an event $A\in\sigma(S_{t^*})$ with $P(A)>0$ such that $\phi_{t^*}^{\theta_n}(S_{t^*})\to0$ on $A$. Assume
\[
\begin{aligned}
\text{(a)}\quad &\inf_n\inf_{\omega\in A}\|\nabla_\theta\phi_{t^*}^{\theta_n}(S_{t^*}(\omega))\|>0,\\
\text{(b)}\quad &\sup_n\sup_{\omega\in A}\phi_s^{\theta_n}(S_s(\omega))\le1-\epsilon_0
                 \quad\text{for every }s<t^*,\\
\text{(c)}\quad &P(\{R_{t^*}<\nu\}\cap A)>0.
\end{aligned}
\]
Then
\[
\E[\|\hat\nabla L_{\rm score}(\theta_n)\|^2]\to\infty,
\qquad
\sup_n\E[\|\hat\nabla L(\theta_n)\|^2]<\infty.
\]
\end{theorem}

\begin{remark}[Sigmoid policies]\label{rem:sigmoid}
We implement the proposed algorithm, DIOS, using sigmoid policies. Fix $\nu\in\R$, let $B=1$, and let $(\theta_n)\subset\Theta$. Assume that the logit gradients are uniformly bounded in $n$, $t$, and the state. Then condition (a) of Theorem~\ref{cor:pathwise-vs-score} fails if $\phi_{t^*}^{\theta_n}(S_{t^*})\to0$. Under Assumptions~\ref{asm:reg},~\ref{asm:smooth}, and~\ref{asm:l2}, if $\phi_t^{\theta_n}(S_t)$ converges almost surely to $0$ or $1$ for every $t<T$, then $\hat\nabla L_{\rm score}(\theta_n)\xrightarrow{L^2}0$ and $\hat\nabla L(\theta_n)\xrightarrow{L^2}0$, as proved in \hyperref[app:sigmoid]{Appendix~\ref*{app:proofs}}.
\end{remark}

\section{Distribution-Informed Optimal Stopping}
\label{sec:dios-section}

DIOS (\emph{Distribution-Informed Optimal Stopping}) combines a per-time-step soft penalty with tail-level annealing. DIOS estimates the reward quantile directly from each batch without an auxiliary critic. Section~\ref{sec:convergence} gives a finite-time stationarity guarantee for the fixed-tail-level update.

\subsection{Blindness to success}
\label{sec:remedies}

CVaR-PG~\citep{tamar2015cvar} estimates the gradient of $\CVaR_\beta(R_{\tau_\theta})$ using a trajectory score. On each of $B$ independent state trajectories $i$, let $\tau_\theta^{(i)}$ be the stopping time generated by $\phi^\theta$, and let $\pi_\theta^{(i)}$ denote the action likelihood defined in Section~\ref{sec:pathwise}, evaluated on trajectory $i$. Let $\widehat\nu_\beta$ be the empirical $\beta$-quantile of $R_{\tau_\theta^{(1)}}^{(1)},\ldots,R_{\tau_\theta^{(B)}}^{(B)}$. The CVaR-PG gradient estimator is
\begin{equation}
\label{eq:cvarpg}
\widehat g_{\rm PG}(\theta)
=\frac1{\beta B}\sum_{i=1}^B
\ind\!\left\{R_{\tau_\theta^{(i)}}^{(i)}\le\widehat\nu_\beta\right\}
\left(R_{\tau_\theta^{(i)}}^{(i)}-\widehat\nu_\beta\right)
\nabla_\theta\log\pi_\theta^{(i)}(\tau_\theta^{(i)}).
\end{equation}
The CVaR-PG estimator $\widehat g_{\rm PG}(\theta)$ in~\eqref{eq:cvarpg} has \emph{blindness to success}: trajectories with sampled stopping rewards above $\widehat\nu_\beta$ receive zero weight.

Methods for improving sample efficiency in CVaR policy optimization modify the gradient estimator, the policy, or the sampling procedure. \citet{greenberg2022} pair a soft-risk schedule with a cross-entropy environment sampler that requires controllable environment dynamics. \citet{luo2024mixture} combine a risk-averse policy with a risk-neutral policy trained via Implicit Q-Learning~\citep{kostrikov2022}. \citet{kim2024pcvar} weight individual costs using an auxiliary tail-probability predictor. \citet{mead2025retcap} compute policy gradients from capped trajectory returns and update the cap during training, subject to an externally supplied minimum cap. \citet{luo2026cvarpg} augment CVaR-PG with a quantile-gradient term computed by an auxiliary quantile critic.

\subsection{The DIOS algorithm}
\label{sec:dios-alg}

For a fixed policy $\phi^\theta$ and $\beta\in(0,1)$, the CVaR envelope~\eqref{eq:ru} is maximized at $\nu=\VaR_\beta(R_{\tau_\theta})$. DIOS, presented in Algorithm~\ref{alg:dios}, samples one stopping time on each trajectory according to the current policy. The reward quantile is then estimated from the batch of rewards obtained at the sampled stopping times. Let $\nu_k$ and $\beta_k$ denote the quantile estimate and tail level at iteration $k$. With $\nu_k$ held fixed, the policy update sums the terms $-\beta_k^{-1}(\nu_k-R_t)^+\nabla_\theta\psi_t^\theta$ over all stopping times on each trajectory. Thus a reward $R_t<\nu_k$ can contribute even when the reward at the sampled stopping time exceeds $\nu_k$. Computing $\nu_k$ directly from the batch removes the auxiliary critic.

Tail-level annealing gradually shifts training toward the target lower tail. The fractions $\rho_w$ and $\rho_a$ specify the warmup and annealing durations, respectively. During warmup, $\beta_k=1$ and $\nu_k$ is the largest of the $B$ sampled stopping rewards. Since $\sum_{t=0}^T\nabla_\theta\psi_t^\theta=0$ and $-(\nu_k-R_t)^+=\min\{R_t,\nu_k\}-\nu_k$, the warmup gradient equals the pathwise gradient with rewards $\min\{R_t,\nu_k\}$. Warmup therefore trains on rewards truncated at the batch maximum of sampled stopping rewards. After warmup, $\beta_k$ decreases linearly to $\beta_{\rm tgt}$ and remains at $\beta_{\rm tgt}$ for the rest of training.

\begin{algorithm}[t]
\caption{DIOS for CVaR optimal stopping}
\label{alg:dios}
\begin{algorithmic}[1]
\REQUIRE policy network $\phi^\theta = (\phi_0^\theta, \dots, \phi_{T-1}^\theta)$ with implicit $\phi_T^\theta \equiv 1$ (forced terminal stop); target tail level $\beta_{\rm tgt} \in (0, 1)$; warmup and annealing fractions $\rho_w,\rho_a\ge0$ with $\rho_w+\rho_a\le1$; positive integers $K,B$; learning rate $\eta>0$
\FOR{$k = 0, 1, \dots, K-1$}
    \STATE If $\rho_a>0$, set $\beta_k \leftarrow 1 + (\beta_{\rm tgt} - 1) \cdot \min\!\left\{1,\max\!\left\{0,\dfrac{k-\rho_wK}{\rho_aK}\right\}\right\}$. If $\rho_a=0$, set $\beta_k=1$ for $k<\rho_wK$ and $\beta_k=\beta_{\rm tgt}$ thereafter.
    \STATE Sample $B$ independent trajectories $\{(S_t^{(i)})_{t=0}^T\}_{i=1}^B$ from the environment
    \STATE Compute $\psi_t^{\theta,(i)} \leftarrow \phi_t^\theta(S_t^{(i)}) \prod_{s<t}\!\left(1 - \phi_s^\theta(S_s^{(i)})\right)$ for all $i, t$
    \STATE Independently for each $i$, sample $\tau^{(i)} \sim \mathrm{Categorical}\!\left((\psi_t^{\theta,(i)})_{t=0}^T\right)$, compute $G^{(i)} := R_{\tau^{(i)}}^{(i)}$
    \STATE $m_k \leftarrow \max\{1,\lfloor\beta_kB\rfloor\}$, \quad $\nu_k \leftarrow$ the $m_k$-th smallest value in $\{G^{(i)}\}_{i=1}^B$
    \STATE $\tilde R_t^{(i)} \leftarrow -\frac{1}{\beta_k}\left(\nu_k - R_t^{(i)}\right)^+$ for all $i, t$
    \STATE $\widehat{L}(\theta) \leftarrow -\frac{1}{B}\sum_{i=1}^B \sum_{t=0}^T \mathrm{sg}\!\left[\tilde R_t^{(i)}\right] \cdot \mathrm{sg}\!\left[\psi_t^{\theta,(i)}\right] \cdot \log \psi_t^{\theta,(i)}$
    \STATE $\theta \leftarrow \theta - \eta \nabla_\theta \widehat{L}(\theta)$
\ENDFOR
\end{algorithmic}
\end{algorithm}

\subsection{Finite-time convergence}
\label{sec:convergence}
We establish a finite-time stationarity guarantee for the DIOS update at a fixed tail level and bound the reward-quantile estimation error.

Fix $\beta\in(0,1)$. Define
\[
V_\beta(\theta):=\CVaR_\beta(R_{\tau_\theta}),\qquad
F_\theta(x):=\Prob(R_{\tau_\theta}\le x).
\]
\begin{samepage}
\begin{assumption}[Convergence regime]\label{asm:conv}
There exist a nonempty open convex domain $U\subset\Theta$ and numbers $\nu_->0$, $f_->0$, $\delta_0>0$, and $L_V>0$. Define
\[
\Theta_{\nu_-}:=\{\theta\in U:\nu^*(\theta)\ge\nu_-\}.
\]
\begin{enumerate}[label=\textbf{(\alph*)},leftmargin=2em,itemsep=1pt,topsep=1pt]
\item For every $\theta\in\Theta_{\nu_-}$ and $u\in[0,\delta_0]$,
\[
F_\theta\bigl(\nu^*(\theta)+u\bigr)\ge\beta+f_-u,
\qquad
F_\theta\bigl(\nu^*(\theta)-u\bigr)\le\beta-f_-u.
\]
\item The function $V_\beta$ is differentiable on $U$, and $\nabla V_\beta$ is $L_V$-Lipschitz on $U$.
\end{enumerate}
\end{assumption}
\end{samepage}
For each integer batch size $B\ge1$, fix an integer $m_B$ satisfying
\begin{equation}
\label{eq:rank}
1\le m_B\le B,\qquad |m_B-\beta B|\le1.
\end{equation}
Let $\mathcal H_k$ be the sigma-field generated by $\theta_0$ and all state trajectories and stopping draws before iteration $k$. Conditional on $\mathcal H_k$, sample $B$ independent copies of $(S_t)_{t=0}^T$. For each trajectory $i$, draw $\tau_{\theta_k}^{(i)}$ with probabilities $(\psi_t^{\theta_k,(i)})_{t=0}^T$, independently across $i$ conditional on the state trajectories. Let $\nu_k$ be the $m_B$-th smallest reward among $R_{\tau_{\theta_k}^{(i)}}^{(i)}$, $i=1,\ldots,B$. For a constant step size $\eta>0$, define
\begin{equation}
\label{eq:update}
\widehat g_k(\nu):=-\frac1{\beta B}\sum_{i=1}^B\sum_{t=0}^T
(\nu-R_t^{(i)})^+\left.\nabla_\theta\psi_t^{\theta,(i)}\right|_{\theta=\theta_k},
\qquad
\theta_{k+1}=\theta_k+\eta\widehat g_k(\nu_k).
\end{equation}
Assumption~\ref{asm:conv}(b) gives~\citep{ghadimi2013}
\[
V_\beta(\theta')\ge V_\beta(\theta)+\langle\nabla V_\beta(\theta),\theta'-\theta\rangle
-\tfrac{L_V}{2}\|\theta'-\theta\|^2,\qquad \theta,\theta'\in U.
\]

\begin{theorem}[Finite-time convergence]\label{thm:conv}
Under Assumptions~\ref{asm:smooth},~\ref{asm:l2}, and~\ref{asm:conv}, there exist an integer $B_0\ge1$ and constants $C_\nu>0$ and $C_g\ge0$, depending only on $\beta$, $T$, $M$, $C_\psi$, $f_-$, and $\delta_0$. Let $K,B$ be integers with $K\ge1$ and $B\ge B_0$, and let $0<\eta\le1/L_V$. Suppose $(\theta_k)_{k=0}^K$ satisfies \eqref{eq:update} and $\theta_k\in\Theta_{\nu_-}$ almost surely for $0\le k\le K$. Define
\[
D_0:=\sup_{\vartheta\in\Theta}V_\beta(\vartheta)-\E[V_\beta(\theta_0)].
\]
Then the following bounds hold.
\begin{enumerate}[label=\textbf{(\alph*)},leftmargin=2em,itemsep=1pt,topsep=1pt]
\item The reward-quantile estimates $\nu_k$ satisfy
\begin{equation}
\label{eq:quantile}
\max_{0\le k<K}\E\!\left[(\nu_k-\nu^*(\theta_k))^2\right]\le\frac{C_\nu}{B}.
\end{equation}
\item The average expected squared gradient satisfies
\begin{equation}
\label{eq:stationarity}
\frac1K\sum_{k=0}^{K-1}\E\!\left[\|\nabla V_\beta(\theta_k)\|^2\right]
\le\frac{2D_0}{\eta K}+\frac{C_g}{B}.
\end{equation}
\end{enumerate}
Consequently,
\begin{equation}
\label{eq:conv}
\min_{0\le k<K}\E\!\left[\|\nabla V_\beta(\theta_k)\|^2+
(\nu_k-\nu^*(\theta_k))^2\right]
\le\frac{2D_0}{\eta K}+\frac{C_g+C_\nu}{B}.
\end{equation}
For every $\epsilon>0$, integers $K\ge1$ and $B\ge B_0$ can be chosen so that
\[
\min_{0\le k<K}\E\!\left[\|\nabla V_\beta(\theta_k)\|^2+(\nu_k-\nu^*(\theta_k))^2\right]\le\epsilon,
\]
provided the corresponding iterates satisfy $\theta_k\in\Theta_{\nu_-}$ almost surely. For fixed $\beta$ and $\eta$ and sufficiently small $\epsilon>0$, the choices satisfy $K=O(\epsilon^{-1})$ and $B=O(\epsilon^{-1})$, giving $KB=O(\epsilon^{-2})$ sampled trajectories.
\end{theorem}

\begin{remark}
Theorem~\ref{thm:conv} applies when quantile estimation and the policy update use the same batch. Since $\nu_k$ depends on the sampled trajectories, $\widehat g_k(\nu_k)$ need not be conditionally unbiased for $\nabla V_\beta(\theta_k)$ given $\mathcal H_k$. The proof in Appendix~\ref{app:convproof} establishes
\[
\E\!\left[\|\widehat g_k(\nu_k)-\nabla V_\beta(\theta_k)\|^2\mid\mathcal H_k\right]\le C_g/B
\]
by bounding the sampling error at $\nu^*(\theta_k)$ and the error from replacing $\nu^*(\theta_k)$ with $\nu_k$.
\end{remark}

\section{Empirical validation}
\label{sec:exp}

We validate DIOS on Bermudan max-call options against four CVaR baselines: CVaR-PG~\citep{tamar2015cvar}, MIX~\citep{luo2024mixture}, RET-CAP~\citep{mead2025retcap}, and CVaR-VaR~\citep{luo2026cvarpg}. Appendix~\ref{app:pcvar} reports additional results for PCVaR-PG~\citep{kim2024pcvar}. The setup uses volatility $\sigma = 0.20$, $d \in \{20, 50, 100, 200\}$ underlying assets, initial price $s_0 \in \{90, 100, 110\}$, $10$ candidate stopping times $t=0,\ldots,T$ with $T=9$, and CVaR level $\beta = 0.2$, yielding 12 configurations. An asymmetric-volatility variant reported in Appendix~\ref{app:asym} contributes 12 more, for 24 in total. We estimate CVaR$_{0.2}$ by averaging the lowest $20\%$ of evaluation rewards. Let $\bar\tau\in[0,T]$ denote the mean stopping time on the evaluation set. For each configuration, the top row of Table~\ref{tab:main} reports the mean and standard deviation of CVaR$_{0.2}$ over ten training seeds. The bottom row reports training time in seconds, including periodic evaluation, and the fraction of seeds with $\bar\tau<0.5$, which we classify as collapsed. Bold entries mark the highest mean CVaR and shortest training time. Implementation details are in Appendix~\ref{app:exp-details}.

\begin{table}[t]
\centering
\caption{Symmetric Bermudan max-call ($\sigma=0.20$). RET-CAP CVaR values use reruns with recomputed minimum caps and retuned hyperparameters. RET-CAP times use the earlier minimum caps and hyperparameters.}
\label{tab:main}
\footnotesize
\setlength{\tabcolsep}{2.2pt}
\begin{tabular}{cc l ccccc}
\toprule
$d$ & $s_0$ & & DIOS (ours) & CVaR-VaR & CVaR-PG & MIX & RET-CAP \\
\midrule
 \multirow{6}{*}{$20$} & \multirow{2}{*}{$90$}
 & CVaR\,$\uparrow$ & $\mathbf{18.171}\,\pm\,0.285$ & $18.068\,\pm\,0.385$ & $14.452\,\pm\,7.232$ & $14.331\,\pm\,7.176$ & $18.003\,\pm\,0.290$ \\
 & & wall\,(coll)\,$\downarrow$ & $\mathbf{12.3}$\,(0/10) & $54.7$\,(0/10) & $14.8$\,(2/10) & $22.8$\,(2/10) & $23.5$\,(0/10) \\
\cmidrule(l){2-8}
  & \multirow{2}{*}{$100$}
 & CVaR\,$\uparrow$ & $\mathbf{29.963}\,\pm\,0.329$ & $29.850\,\pm\,0.441$ & $29.872\,\pm\,0.535$ & $28.004\,\pm\,5.385$ & $29.677\,\pm\,0.295$ \\
 & & wall\,(coll)\,$\downarrow$ & $\mathbf{12.7}$\,(0/10) & $51.0$\,(0/10) & $15.2$\,(0/10) & $24.8$\,(0/10) & $25.8$\,(0/10) \\
\cmidrule(l){2-8}
  & \multirow{2}{*}{$110$}
 & CVaR\,$\uparrow$ & $\mathbf{41.966}\,\pm\,0.339$ & $41.893\,\pm\,0.587$ & $32.222\,\pm\,14.557$ & $25.787\,\pm\,15.793$ & $41.616\,\pm\,0.366$ \\
 & & wall\,(coll)\,$\downarrow$ & $\mathbf{13.0}$\,(0/10) & $51.6$\,(0/10) & $14.3$\,(3/10) & $23.3$\,(5/10) & $26.0$\,(0/10) \\
\midrule
 \multirow{6}{*}{$50$} & \multirow{2}{*}{$90$}
 & CVaR\,$\uparrow$ & $\mathbf{33.220}\,\pm\,0.126$ & $33.213\,\pm\,0.175$ & $33.121\,\pm\,0.144$ & $33.146\,\pm\,0.173$ & $32.810\,\pm\,0.225$ \\
 & & wall\,(coll)\,$\downarrow$ & $\mathbf{11.2}$\,(0/10) & $48.8$\,(0/10) & $15.3$\,(0/10) & $23.9$\,(0/10) & $26.0$\,(0/10) \\
\cmidrule(l){2-8}
  & \multirow{2}{*}{$100$}
 & CVaR\,$\uparrow$ & $\mathbf{46.802}\,\pm\,0.153$ & $46.715\,\pm\,0.187$ & $46.614\,\pm\,0.272$ & $46.588\,\pm\,0.289$ & $45.913\,\pm\,0.278$ \\
 & & wall\,(coll)\,$\downarrow$ & $\mathbf{11.4}$\,(0/10) & $58.2$\,(0/10) & $15.0$\,(0/10) & $22.0$\,(0/10) & $25.1$\,(0/10) \\
\cmidrule(l){2-8}
  & \multirow{2}{*}{$110$}
 & CVaR\,$\uparrow$ & $\mathbf{60.392}\,\pm\,0.171$ & $60.241\,\pm\,0.311$ & $60.157\,\pm\,0.275$ & $60.257\,\pm\,0.234$ & $59.115\,\pm\,0.379$ \\
 & & wall\,(coll)\,$\downarrow$ & $\mathbf{11.6}$\,(0/10) & $59.6$\,(0/10) & $15.2$\,(0/10) & $25.0$\,(0/10) & $25.8$\,(0/10) \\
\midrule
 \multirow{6}{*}{$100$} & \multirow{2}{*}{$90$}
 & CVaR\,$\uparrow$ & $\mathbf{45.001}\,\pm\,0.167$ & $44.859\,\pm\,0.195$ & $40.406\,\pm\,13.470$ & $40.408\,\pm\,13.471$ & $44.323\,\pm\,0.359$ \\
 & & wall\,(coll)\,$\downarrow$ & $\mathbf{11.9}$\,(0/10) & $52.6$\,(0/10) & $15.1$\,(1/10) & $22.5$\,(1/10) & $25.4$\,(0/10) \\
\cmidrule(l){2-8}
  & \multirow{2}{*}{$100$}
 & CVaR\,$\uparrow$ & $\mathbf{59.838}\,\pm\,0.160$ & $59.805\,\pm\,0.148$ & $47.853\,\pm\,23.927$ & $47.771\,\pm\,23.886$ & $58.956\,\pm\,0.390$ \\
 & & wall\,(coll)\,$\downarrow$ & $\mathbf{11.9}$\,(0/10) & $49.9$\,(0/10) & $14.9$\,(2/10) & $22.9$\,(2/10) & $25.4$\,(0/10) \\
\cmidrule(l){2-8}
  & \multirow{2}{*}{$110$}
 & CVaR\,$\uparrow$ & $\mathbf{74.668}\,\pm\,0.161$ & $74.543\,\pm\,0.202$ & $61.470\,\pm\,25.738$ & $61.580\,\pm\,25.791$ & $73.084\,\pm\,0.600$ \\
 & & wall\,(coll)\,$\downarrow$ & $\mathbf{12.0}$\,(0/10) & $59.5$\,(0/10) & $15.6$\,(2/10) & $22.5$\,(2/10) & $25.3$\,(0/10) \\
\midrule
 \multirow{6}{*}{$200$} & \multirow{2}{*}{$90$}
 & CVaR\,$\uparrow$ & $\mathbf{57.142}\,\pm\,0.108$ & $57.110\,\pm\,0.117$ & $51.275\,\pm\,17.093$ & $51.302\,\pm\,17.102$ & $56.105\,\pm\,0.423$ \\
 & & wall\,(coll)\,$\downarrow$ & $\mathbf{11.8}$\,(0/10) & $59.0$\,(0/10) & $14.3$\,(1/10) & $25.0$\,(1/10) & $26.8$\,(0/10) \\
\cmidrule(l){2-8}
  & \multirow{2}{*}{$100$}
 & CVaR\,$\uparrow$ & $\mathbf{73.259}\,\pm\,0.074$ & $73.233\,\pm\,0.126$ & $65.788\,\pm\,21.930$ & $58.529\,\pm\,29.264$ & $71.396\,\pm\,0.410$ \\
 & & wall\,(coll)\,$\downarrow$ & $\mathbf{11.8}$\,(0/10) & $58.6$\,(0/10) & $15.7$\,(1/10) & $22.2$\,(2/10) & $25.4$\,(0/10) \\
\cmidrule(l){2-8}
  & \multirow{2}{*}{$110$}
 & CVaR\,$\uparrow$ & $\mathbf{89.378}\,\pm\,0.100$ & $89.375\,\pm\,0.099$ & $81.387\,\pm\,23.796$ & $81.287\,\pm\,23.764$ & $86.622\,\pm\,0.795$ \\
 & & wall\,(coll)\,$\downarrow$ & $\mathbf{11.7}$\,(0/10) & $57.0$\,(0/10) & $15.8$\,(1/10) & $23.3$\,(1/10) & $25.3$\,(0/10) \\
\bottomrule
\end{tabular}
\end{table}

DIOS achieves the highest mean CVaR and the shortest recorded training time in all 24 symmetric and asymmetric configurations. CVaR-PG and MIX have large standard deviations and collapsed seeds, with collapse rates of $3/10$ and $5/10$, respectively, at symmetric $(d,s_0)=(20,110)$. CVaR-VaR takes roughly $4$--$5$ times as long as DIOS, while MIX and RET-CAP take roughly twice as long. Appendix~\ref{app:cvarvar-critic} compares CVaR-VaR with a variant that replaces the quantile critic by a batch mean and changes the tradeoff-weight schedule. For RET-CAP, we compute the minimum cap $q^*$ by averaging empirical $\beta$-quantiles of DIOS evaluation rewards over ten seeds. Appendix~\ref{app:retcap-anchor} specifies the calculation and Appendix~\ref{app:asym} reports the asymmetric results.

The mean advantage over CVaR-VaR is small in some configurations. At symmetric $(d,s_0)=(200,110)$, mean CVaR is $89.378$ for DIOS and $89.375$ for CVaR-VaR. Two-sided paired $t$-tests over ten training seeds give $p<0.05$ against CVaR-VaR in 5 of the 24 configurations. The $p$-values are not adjusted for testing multiple configurations. Appendix~\ref{app:paired} gives the paired mean differences and $p$-values.

\paragraph{Additional experiments.}
Appendix~\ref{app:smallbeta} compares raw batch quantile estimates with exponential moving average (EMA) estimates within DIOS at $\beta=0.1$ and $0.05$. DIOS with EMA has the highest mean CVaR among the tested methods in all 12 configurations at each tail level. Paired $t$-tests give $p<0.05$ against CVaR-VaR in 6/12 configurations at $\beta=0.1$ and 12/12 at $\beta=0.05$, with no adjustment for multiple comparisons. Appendix~\ref{app:house} reports house-selling experiments, where DIOS has the highest mean CVaR in 23/24 configurations. Appendix~\ref{app:width} compares CVaR, parameter counts, and training times for DIOS encoder widths 32--256.

\FloatBarrier

\section{Conclusion}
\label{sec:conclusion}

We developed a framework and algorithm for CVaR deep optimal stopping. The properties (S1)--(S3) yield Markov optimality without state augmentation and a pathwise gradient formula. Theorem~\ref{cor:pathwise-vs-score} contrasts bounded pathwise second moments with score-function divergence under the stated conditions. DIOS combines a per-time-step soft penalty, tail-level annealing, and direct batch quantile estimation. We establish a finite-time stationarity guarantee at a fixed tail level. DIOS achieves the highest mean CVaR and the shortest recorded training time across the 24 max-call configurations.

\paragraph{Limitations and future work.}
Theorem~\ref{thm:conv} bounds expected squared gradients and squared quantile errors at fixed $\beta$ under Assumptions~\ref{asm:smooth},~\ref{asm:l2}, and~\ref{asm:conv}, with $\theta_k\in\Theta_{\nu_-}$ throughout the analyzed updates. The theorem does not establish parameter convergence or global optimality. The score-function and pathwise estimators can each vanish under the sigmoid conditions in \hyperref[app:sigmoid]{Appendix~\ref*{app:proofs}}. Appendix~\ref{app:exp-verif} explains why ReLU, unbounded inputs, training-mode BatchNorm, and reward atoms prevent automatic verification of the assumptions. Extending the convergence analysis to Adam, clipping, training-mode BatchNorm, annealing, and EMA remains open. At $\beta=1$, unbounded rewards can preclude a finite maximizing $\nu$ in~\eqref{eq:ru}.

The annealing ablation also changes the quantile rank, and the critic ablation also changes the tradeoff-weight schedule. Recorded training times come from different sessions, sometimes on shared GPUs, and exclude training and evaluation of policies used to set RET-CAP's minimum cap. RET-CAP CVaR values and times use different settings, as specified in Table~\ref{tab:main}. The experiments cover max-call options, house selling, and encoder widths 32--256. Further work includes larger architectures, action-dependent dynamics, and convergence for other optimized certainty equivalent objectives. Medical and control applications require validation of tail estimates to limit human and financial costs.

\FloatBarrier

\clearpage
\section*{Acknowledgments}
This work was supported by the Korea government (MSIT) through the NRF (RS-2023-00219980 \& RS-2026-25491270); and by the Institute of Information \& Communications Technology Planning \& Evaluation (IITP) grant funded by the Korea government (MSIT) (RS-2021-II211343, Artificial Intelligence Graduate School Program, Seoul National University).


\appendix

\section{Assumptions in the experimental domains}
\label{app:exp-verif}

We verify the payoff-moment assumptions and examine the additional policy and convergence conditions in the experimental domains.

\subsection{Model and payoff moments}
\paragraph{Max-call rewards.}
Let $S_t^j$ denote the price of asset $j\in\{1,\ldots,d\}$, with $S_0^j=s_0$. Asset prices follow independent geometric Brownian motions under the risk-neutral measure, with drift $r-\delta$ and volatility $\sigma_j$. Prices are sampled at $t=0,\ldots,T=9$ with time step $\Delta t=1/3$ year. Let $K_{\rm strike}=100$ and define
\[
R_t:=e^{-rt\Delta t}\left(\max_{1\le j\le d}S_t^j-K_{\rm strike}\right)^+.
\]
The experiments use $r=0.05$, $\delta=0.10$, $d\in\{20,50,100,200\}$, and $s_0\in\{90,100,110\}$. The symmetric domain has $\sigma_j=0.20$. The asymmetric domain has $\sigma_j=0.10+j/(2d)$. Every asset price has finite moments of all orders. Since
\[
0\le R_t\le\sum_{j=1}^d S_t^j,
\qquad R_t^2\le\sum_{j=1}^d(S_t^j)^2,
\]
Assumptions~\ref{asm:reg} and~\ref{asm:l2} hold in both domains, with $M=\max_{0\le t\le T}\|R_t\|_{L^2(P)}<\infty$.

\paragraph{House-selling rewards.}
Let $X_j>0$ be the offer at decision $j\in\{1,\ldots,9\}$ and let $c>0$ be the cost per decision. Accepting at $j$ gives $X_j-cj$, which can be negative. Adding $9c$ to the reward at every decision makes the rewards nonnegative. For each fixed policy, CVaR and VaR increase by $9c$. Replacing $R_t$ by $R_t+9c$ and $\nu$ by $\nu+9c$ leaves $(\nu-R_t)^+$ and the pathwise gradient unchanged. The tables retain the original reward units. Lognormal offers and the tested Pareto shapes $a>2$ have finite second moments.

\subsection{Policy regularity}
\paragraph{Sigmoid output.}
For a finite logit $z_t^\theta(S_t)$, the probability $\phi_t^\theta(S_t)=(1+e^{-z_t^\theta(S_t)})^{-1}$ lies in $(0,1)$ for $t<T$. Forced terminal stopping gives $\phi_T^\theta(S_T)=1$ and $\nabla_\theta\phi_T^\theta(S_T)=0$. A sigmoid output does not by itself imply the uniform derivative bound in Assumption~\ref{asm:smooth}. For example, $\phi^\theta(s)=(1+e^{-\theta s})^{-1}$ has derivative $s/4$ at $\theta=0$, which is unbounded as $s\to\infty$.

\paragraph{Derivatives of stopping probabilities.}
Suppose $\|\nabla_\theta\phi_t^\theta(S_t)\|\le G$ uniformly in $t$ and $\theta$, as required by Assumption~\ref{asm:smooth}. The derivative formula in Remark~\ref{rem:boundary} contains probability factors in $[0,1]$. The triangle inequality therefore gives
\begin{equation}
\label{eq:Cpsi-loose}
\|\nabla_\theta\psi_t^\theta\|\le G(1+t)\le G(T+1).
\end{equation}
Thus $C_\psi=G(T+1)$ is a valid choice in Assumption~\ref{asm:smooth}.
For a bound independent of $T$, suppose additionally that $\phi_s^\theta(S_s)\ge\phi_-$ for every $s<T$ and $\theta\in\Theta$ on an event of probability one, where $\phi_-\in(0,1)$. For $1\le t\le T$,
\[
\|\nabla_\theta\psi_t^\theta\|
\le G(1-\phi_-)^t+tG(1-\phi_-)^{t-1}.
\]
At $t=0$, $\|\nabla_\theta\psi_0^\theta\|\le G$. Consequently, a valid choice is
\[
C_\psi=G\max\left\{1,\sup_{m\in\{1,2,\ldots\}}
\bigl[(1-\phi_-)^m+m(1-\phi_-)^{m-1}\bigr]\right\}<\infty.
\]
Finiteness follows from $m(1-\phi_-)^{m-1}\to0$. For fixed $G$ and $\phi_-$, the bound $C_\psi$ is independent of $T$.

\paragraph{Training implementation.}
The network in Appendix~\ref{app:exp-details} uses ReLU and BatchNorm. ReLU is not differentiable at zero, and the lognormal inputs are unbounded. During training, BatchNorm statistics use all $B(T+1)$ state-time pairs, so a probability at time $t$ can depend on future states and other trajectories. The training computation therefore differs from the pointwise policy $\phi_t^\theta(S_t)$ assumed in the analysis. Evaluation uses fixed running statistics. The architecture alone does not establish Assumption~\ref{asm:smooth}.

\subsection{Convergence assumptions and training diagnostics}
\paragraph{Terminal-stop benchmark.}
For nonnegative $R_T$, $\CVaR_\beta(R_T)>0$ if and only if $P(R_T=0)<\beta$. Let $\Phi$ denote the standard normal CDF. In the symmetric domain, $T\Delta t=3$ and $r-\delta-\sigma^2/2=-0.07$, giving
\[
P(R_T=0)=\Phi\!\left(\frac{\log(K_{\rm strike}/s_0)+0.21}{0.2\sqrt3}\right)^d.
\]
The largest probability on the tested grid occurs at $(d,s_0)=(20,90)$ and is below $0.0183$. In the asymmetric domain,
\[
P(R_T=0)=\prod_{j=1}^d\Phi\!\left(
\frac{\log(K_{\rm strike}/s_0)-3(r-\delta-\sigma_j^2/2)}{\sigma_j\sqrt3}\right).
\]
The largest probability again occurs at $(d,s_0)=(20,90)$ and is below $0.00861$. Hence stopping at $T$ gives positive CVaR at every tested target tail level in both domains.

\paragraph{Conditions along the iterates.}
The terminal-stop calculation establishes a positive-CVaR benchmark. Theorem~\ref{thm:conv} additionally requires the CDF bounds in Assumption~\ref{asm:conv}(a), a differentiable $V_\beta$ with an $L_V$-Lipschitz gradient on $U$, and $\nu^*(\theta_k)\ge\nu_-$ throughout the analyzed updates. Max-call rewards have an atom at zero, so the CDF bounds cannot be inferred from the continuous asset-price distributions. During the fixed-tail phase, the recorded batch estimate $\nu_k$ is zero in 6 of 180,000 symmetric-domain updates and 2 of 180,000 asymmetric-domain updates. Counts of zero $\nu_k$ do not verify $\nu^*(\theta_k)\ge\nu_-$.

\section{Proofs}
\label{app:proofs}

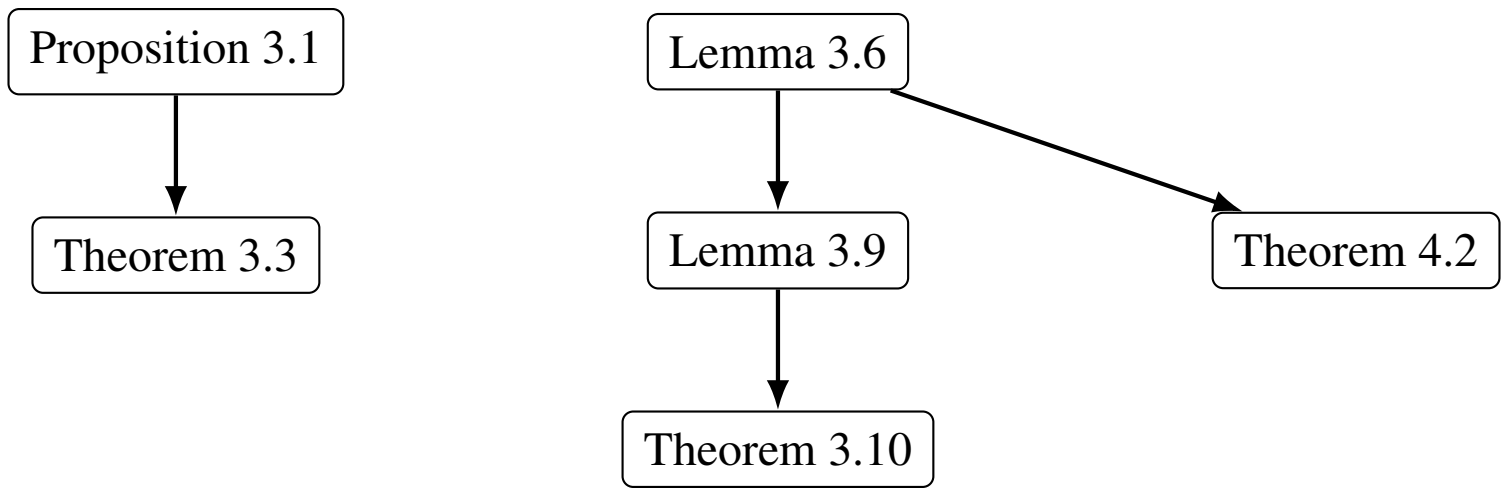
\begin{figure}[ht]
\centering
\begin{tikzpicture}[
  node distance=8mm and 14mm,
  every node/.style={draw, rounded corners=2pt, inner sep=4pt, font=\small},
  arr/.style={-{Latex[length=2mm]}, thick}
]
\node (env) {Proposition~\ref{thm:envelope}};
\node[below=of env] (markov) {Theorem~\ref{thm:markov}};
\node[right=20mm of env] (gradient) {Lemma~\ref{thm:pathwise}};
\node[below=of gradient] (moment) {Lemma~\ref{thm:variance}};
\node[below=of moment] (compare) {Theorem~\ref{cor:pathwise-vs-score}};
\node[right=20mm of moment] (conv) {Theorem~\ref{thm:conv}};
\draw[arr] (env) -- (markov);
\draw[arr] (gradient) -- (moment);
\draw[arr] (moment) -- (compare);
\draw[arr] (gradient) -- (conv);
\end{tikzpicture}
\caption{Proof dependencies. An arrow points from a result to a proof that uses it.}
\label{fig:dep}
\end{figure}

Figure~\ref{fig:dep} shows which preceding results each proof uses.

\subsection*{Proof of Proposition~\ref{thm:envelope}}

\begin{proof}

For every $X \in L^1(P)$, the representation~\eqref{eq:ru} gives
\[
\CVaR_\beta(X) \;=\; \sup_{\nu \in \R} F(\nu, X) \quad \text{with} \quad F(\nu, X) \;:=\; \nu + \tfrac{1}{\beta}\, \E\!\left[-(\nu - X)^+\right],
\]
where the supremum is attained at $\nu = \VaR_\beta(X)$. Assumption~\ref{asm:reg} gives $R_\tau \in L^1(P)$ for every $\tau \in \cT$, so we may set $X = R_\tau$.

By the commutativity of iterated suprema,
\[
V^*_\beta \;=\; \sup_{\tau \in \cT} \sup_{\nu \in \R} F(\nu, R_\tau) \;=\; \sup_{\nu \in \R} \Phi(\nu),
\]
where
\[
\Phi(\nu) \;:=\; \sup_{\tau \in \cT} f_\tau(\nu), \qquad f_\tau(\nu) \;:=\; \nu + \tfrac{1}{\beta}\, \E[-(\nu - R_\tau)^+].
\]
Equation~\eqref{eq:envelope} follows.

For each $x \in \R$, $\nu \mapsto -(\nu - x)^+$ is $1$-Lipschitz:
\[
\bigl|-(\nu_1 - x)^+ + (\nu_2 - x)^+\bigr| \;\le\; |\nu_1 - \nu_2|.
\]
Hence $\nu \mapsto \E[-(\nu - R_\tau)^+]$ is $1$-Lipschitz for every $\tau \in \cT$, and $f_\tau$ is $(1 + 1/\beta)$-Lipschitz. Since $-(\nu - x)^+ \le 0$, $f_\tau(\nu) \le \nu$, so $\Phi(\nu) \le \nu < \infty$ for every $\nu \in \R$. A pointwise-finite supremum of $L$-Lipschitz functions is again $L$-Lipschitz, so $\Phi$ is $(1 + 1/\beta)$-Lipschitz.

For $\nu \le 0$, $R_\tau \ge 0$ (since each $g_t$ takes values in $\R_+$) gives $(\nu - R_\tau)^+ = 0$, hence $\E[-(\nu - R_\tau)^+] = 0$ for every $\tau \in \cT$ and
\[
\Phi(\nu) \;=\; \nu \qquad (\nu \le 0).
\]
In particular $\Phi(0) = 0$, and $\Phi(\nu) \to -\infty$ as $\nu \to -\infty$.

The bound $-(\nu - R_\tau)^+ \le -(\nu - R_\tau)$ is pointwise, so
\[
\Phi(\nu) \;\le\; \nu + \tfrac{1}{\beta}\!\left(-\nu + \sup_{\tau \in \cT} \E[R_\tau]\right) \;=\; (1 - 1/\beta)\,\nu \;+\; \tfrac{1}{\beta}\, \sup_{\tau \in \cT} \E[R_\tau].
\]
Since $R_t \ge 0$, $R_\tau = \sum_{t=0}^T R_t \mathbf{1}\{\tau = t\} \le \sum_{t=0}^T R_t$, hence
\[
\sup_{\tau \in \cT} \E[R_\tau] \;\le\; \sum_{t = 0}^T \E[R_t] \;=:\; C \;<\; \infty
\]
by Assumption~\ref{asm:reg}. Hence $\Phi(\nu) \le (1 - 1/\beta)\,\nu + C/\beta$, and since $1 - 1/\beta < 0$, $\Phi(\nu) \to -\infty$ as $\nu \to +\infty$.

$\Phi$ is continuous and $\Phi(\nu) \to -\infty$ as $|\nu| \to \infty$, so there is $K > 0$ with $\Phi(\nu) < 0$ for $|\nu| > K$. Since $\Phi(0) = 0$, $\sup_\nu \Phi = \sup_{[-K, K]} \Phi$, attained by Weierstrass. Fix $\bar\nu \in [-K, K]$ with $\Phi(\bar\nu) = V^*_\beta$.

Fix $\nu \in \R$. The process $h^\nu_t := -(\nu - R_t)^+$ is $(\cF_t)$-adapted with $|h^\nu_t| \le |\nu| + R_t \in L^1(P)$. Define the Snell envelope
\[
V^\nu_T \;:=\; h^\nu_T, \qquad V^\nu_t \;:=\; \max\bigl(h^\nu_t,\ \E[V^\nu_{t+1} \mid \cF_t]\bigr) \quad (t < T).
\]
By~\citet{shiryaev2008optimal},
\[
\sup_{\tau \in \cT} \E[h^\nu_\tau] \;=\; \E[V^\nu_0],
\]
where the supremum is attained at $\tau^\nu := \inf\{t \ge 0 : V^\nu_t = h^\nu_t\}$. The boundary condition $V^\nu_T = h^\nu_T$ forces $\tau^\nu \le T$, so $\tau^\nu \in \cT$. Since $f_\tau(\nu) = \nu + \tfrac{1}{\beta}\,\E[h^\nu_\tau]$, $\Phi(\nu) = f_{\tau^\nu}(\nu)$.

Set $\tau^* := \tau^{\bar\nu}$. Then
\[
V^*_\beta \;=\; \Phi(\bar\nu) \;=\; f_{\tau^*}(\bar\nu) \;\le\; \sup_\nu f_{\tau^*}(\nu) \;=\; \CVaR_\beta(R_{\tau^*}) \;\le\; V^*_\beta.
\]
By~\eqref{eq:ru}, $\sup_\nu f_{\tau^*}(\nu)$ is attained at $\nu = \VaR_\beta(R_{\tau^*})$, so
\[
f_{\tau^*}\bigl(\VaR_\beta(R_{\tau^*})\bigr) \;=\; V^*_\beta.
\]
Hence $\Phi(\VaR_\beta(R_{\tau^*}))\ge V^*_\beta$. The definition $V^*_\beta=\sup_{\nu\in\R}\Phi(\nu)$ gives the reverse inequality. Setting $\nu^* := \VaR_\beta(R_{\tau^*})$, the pair $(\nu^*, \tau^*)$ jointly realizes the supremum in~\eqref{eq:envelope}.
\end{proof}

\subsection*{Proof of Theorem~\ref{thm:markov}}

\begin{proof}

By Proposition~\ref{thm:envelope}, there exists $\nu^* \in \R$ attaining the outer supremum in~\eqref{eq:envelope}. Define $h^{\nu^*}_t := -(\nu^* - R_t)^+$ and the Snell envelope
\[
V^{\nu^*}_T \;:=\; h^{\nu^*}_T, \qquad V^{\nu^*}_t \;:=\; \max\bigl(h^{\nu^*}_t,\ \E[V^{\nu^*}_{t+1} \mid \cF_t]\bigr) \quad (t < T).
\]
By~\citet{shiryaev2008optimal}, $\tau^{\nu^*} := \inf\{t \ge 0 : V^{\nu^*}_t = h^{\nu^*}_t\}$ attains $\sup_{\tau \in \cT} \E[h^{\nu^*}_\tau]$.

It suffices to find Borel sets $D_0, \dots, D_T \subseteq \R^d$ with $D_T = \R^d$ and
\[
\tau^{\nu^*} \;=\; \inf\{t \ge 0 : S_t \in D_t\} \quad P\text{-a.s.}
\]
Define the Borel function $\tilde h^{\nu^*}_t:\R^d\to\R$ by $\tilde h^{\nu^*}_t(s):=-(\nu^*-g_t(s))^+$. By (S1), $h^{\nu^*}_t=\tilde h^{\nu^*}_t(S_t)$.

Backward induction gives Borel functions $v^{\nu^*}_t:\R^d\to\R$ with $v^{\nu^*}_t(S_t)\in L^1(P)$ satisfying
\[
V^{\nu^*}_t \;=\; v^{\nu^*}_t(S_t) \quad P\text{-a.s.}
\]
At $t=T$, define $v^{\nu^*}_T:=\tilde h^{\nu^*}_T$. Then $V^{\nu^*}_T=v^{\nu^*}_T(S_T)$ and $v^{\nu^*}_T(S_T)\in L^1(P)$ by Assumption~\ref{asm:reg}.

Suppose $V^{\nu^*}_{t+1} = v^{\nu^*}_{t+1}(S_{t+1})$ $P$-a.s.\ for some Borel $v^{\nu^*}_{t+1}$ with $v^{\nu^*}_{t+1}(S_{t+1}) \in L^1(P)$. Since $(S_t)_{t=0}^T$ is Markov,
\[
\E[V^{\nu^*}_{t+1} \mid \cF_t] \;=\; \E[v^{\nu^*}_{t+1}(S_{t+1}) \mid S_t].
\]
By Doob-Dynkin, there is a Borel $q^{\nu^*}_t : \R^d \to \R$ with $\E[v^{\nu^*}_{t+1}(S_{t+1}) \mid S_t] = q^{\nu^*}_t(S_t)$ $P$-a.s. Define
\[
v^{\nu^*}_t \;:=\; \max\bigl(\tilde h^{\nu^*}_t,\ q^{\nu^*}_t\bigr),
\]
which is Borel as the max of two Borel functions. Then $V^{\nu^*}_t = v^{\nu^*}_t(S_t)$ $P$-a.s. Integrability follows from
\[
|v^{\nu^*}_t(S_t)| \;\le\; |\tilde h^{\nu^*}_t(S_t)| + |q^{\nu^*}_t(S_t)|
\]
combined with Assumption~\ref{asm:reg} and the conditional Jensen inequality.

Define
\[
D_t \;:=\; \{s \in \R^d : v^{\nu^*}_t(s) = \tilde h^{\nu^*}_t(s)\},
\]
which is Borel as the preimage of $\{0\}$ under the Borel map $v^{\nu^*}_t - \tilde h^{\nu^*}_t$. Then $V^{\nu^*}_t = h^{\nu^*}_t \iff S_t \in D_t$ $P$-a.s., so
\[
\tau^{\nu^*} \;=\; \inf\{t \ge 0 : S_t \in D_t\} \quad P\text{-a.s.}
\]
Since $v^{\nu^*}_T = \tilde h^{\nu^*}_T$, $D_T = \R^d$. Define $\tau^*:=\tau^{\nu^*}$. The stopping time $\tau^*$ is Markov by Definition~\ref{def:markov-stopping}. Since $\nu^*$ attains the outer supremum in~\eqref{eq:envelope}, $V^*_\beta = \nu^* + \tfrac{1}{\beta}\sup_\tau \E[h^{\nu^*}_\tau]$. Since $\tau^*=\tau^{\nu^*}$ maximizes $\E[h^{\nu^*}_\tau]$ over $\tau\in\cT$, $V^*_\beta = \nu^* + \tfrac{1}{\beta}\E[-(\nu^* - R_{\tau^*})^+]$.
\end{proof}

\subsection*{Proof of Lemma~\ref{thm:pathwise}}

\begin{proof}

Since $\psi_t^\theta\, \nabla_\theta \log \psi_t^\theta = \nabla_\theta \psi_t^\theta$,
\[
-\nabla_\theta L(\theta) \;=\; \E\!\left[\sum_{t=0}^T \tilde R^\nu_t\, \nabla_\theta \psi_t^\theta\right].
\]

By (S2), the joint distribution of $(S_0,\ldots,S_T)$ does not depend on $\theta$. At fixed $\nu$, the factors $\tilde R_t^\nu$ also do not depend on $\theta$. Since $\|\nabla_\theta \psi_t^\theta \cdot \tilde R^\nu_t\| \le \tfrac{C_\psi}{\beta}(|\nu| + R_t) \in L^1(P)$ by Assumptions~\ref{asm:reg} and~\ref{asm:smooth}, differentiation under the expectation gives
\[
\nabla_\theta J^{(\nu)}(\theta) \;=\; \E\!\left[\sum_{t=0}^T \tilde R^\nu_t\, \nabla_\theta \psi_t^\theta\right] \;=\; -\nabla_\theta L(\theta).
\]
\end{proof}

\subsection*{Proof of Lemma~\ref{thm:variance}}

\begin{proof}

Fix $\theta\in\Theta$ and $\nu\in\R$. For each sampled trajectory $i=1,\ldots,B$, define
\[
\hat g^{(i)}(\theta) \;:=\; -\sum_{t=0}^T \tilde R^{\nu,(i)}_t \cdot \nabla_\theta \psi_t^{\theta,(i)}.
\]

Define $\bar a^{(i)} := \tfrac{1}{T+1} \sum_{t=0}^T \tilde R^{\nu,(i)}_t$. By (S3), $\sum_{t=0}^T \nabla_\theta \psi_t^{\theta,(i)} = 0$ $P$-a.s., so
\[
\hat g^{(i)}(\theta) \;=\; -\sum_{t=0}^T \bigl(\tilde R^{\nu,(i)}_t - \bar a^{(i)}\bigr) \cdot \nabla_\theta \psi_t^{\theta,(i)}.
\]

By Cauchy--Schwarz,
\[
\|\hat g^{(i)}\|^2 \;\le\; \Bigl(\sum_{t=0}^T (\tilde R^{\nu,(i)}_t - \bar a^{(i)})^2\Bigr) \cdot \Bigl(\sum_{t=0}^T \|\nabla_\theta \psi_t^{\theta,(i)}\|^2\Bigr).
\]

By Assumption~\ref{asm:smooth}, $\|\nabla_\theta \psi_t^{\theta,(i)}\| \le C_\psi$ $P$-a.s., hence
\[
\sum_{t=0}^T \|\nabla_\theta \psi_t^{\theta,(i)}\|^2 \;\le\; (T+1) C_\psi^2.
\]

Since $\bar a^{(i)}$ minimizes $\sum_{t=0}^T(\tilde R_t^{\nu,(i)}-a)^2$ over $a\in\R$, comparison with $a=0$ gives
\[
\E\!\left[\sum_{t=0}^T(\tilde R_t^{\nu,(i)}-\bar a^{(i)})^2\right]
\le\sum_{t=0}^T\E[(\tilde R_t^{\nu,(i)})^2]\le(T+1)\tilde M^2.
\]
Thus $\E[\|\hat g^{(i)}\|^2]\le C_\psi^2(T+1)^2\tilde M^2$.
By convexity, $\E[\|\hat\nabla L(\theta)\|^2]\le B^{-1}\sum_{i=1}^B\E[\|\hat g^{(i)}(\theta)\|^2]\le C_\psi^2(T+1)^2\tilde M^2$, proving part (a).

$\hat\nabla L(\theta) = (1/B) \sum_{i=1}^B \hat g^{(i)}$ is the mean of $B$ i.i.d.\ samples with $\E[\hat g^{(i)}] = \nabla L(\theta)$ by Lemma~\ref{thm:pathwise}, so
\[
\E\bigl[\|\hat\nabla L(\theta) - \nabla L(\theta)\|^2\bigr] \;=\; \tfrac{1}{B}\,\E[\|\hat g^{(1)}-\E[\hat g^{(1)}]\|^2] \;\le\; \tfrac{1}{B}\, \E\bigl[\|\hat g^{(1)}\|^2\bigr] \;\le\; \frac{C_\psi^2(T+1)^2\tilde M^2}{B}.
\]
The variance bound proves part (b). For part (c), if $M=0$, then $R_t=0$ $P$-a.s.\ for every $t$. Thus $\nu^*(\theta)=0$, $\hat\nabla L(\theta)=\nabla L(\theta)=0$ at $\nu=\nu^*(\theta)$, and~\eqref{eq:variance-bound} holds with $c_\beta=0$. Suppose $M>0$. At $\nu=\nu^*(\theta)$, $\tilde M^2\le c_\beta M^2/\beta^2$ and $(T+1)^2\le4T^2$ for $T\ge1$, which gives~\eqref{eq:variance-bound}. For $T=0$, $\psi_0=1$ and $\hat\nabla L(\theta)=\nabla L(\theta)=0$.
For $M>0$, the lower-quantile property gives
\[
(1-\beta)\nu^*(\theta)^2\le\E[R_{\tau_\theta}^2]\le(T+1)M^2.
\]
Since $R_t\ge0$, $((\nu^*(\theta)-R_t)^+)^2\le\nu^*(\theta)^2$, which gives $c_\beta\le(T+1)/(1-\beta)$.
\end{proof}

\subsection*{Proof of Theorem~\ref{cor:pathwise-vs-score}}

\begin{proof}

Define $\tilde M:=\max_t\|\tilde R_t^\nu\|_{L^2(P)}$.

On $\{\tau=u\}$, the sampled action likelihood equals $\psi_u^\theta$. For $u<T$,
\[
\nabla_\theta\log\psi_u^\theta \;=\; \frac{\nabla\phi_u^\theta(S_u)}{\phi_u^\theta(S_u)} \;-\; \sum_{s < u} \frac{\nabla\phi_s^\theta(S_s)}{1 - \phi_s^\theta(S_s)}.
\]
At $u=T$,
\[
\nabla_\theta\log\psi_T^\theta
=-\sum_{s<T}\frac{\nabla\phi_s^\theta(S_s)}{1-\phi_s^\theta(S_s)}.
\]
For each $u$, $\nabla_\theta\log\psi_u^\theta$ is $\cF_u$-measurable.

Since $\E[\mathbf{1}\{\tau = u\} \mid \cF_T] = \phi_u^\theta(S_u) \prod_{s<u}(1-\phi_s^\theta(S_s)) = \psi_u^\theta$,
\[
\E\bigl[\|\hat\nabla L_{\rm score}(\theta_n)\|^2\bigr] \;=\; \sum_{u=0}^T \E\!\left[(\tilde R^\nu_u)^2\, \|\nabla_\theta\log\psi_u^{\theta_n}\|^2\, \psi_u^{\theta_n}\right].
\]

By condition (c) and continuity of probability under increasing unions, choose $\eta>0$ such that $P(\{R_{t^*}\le\nu-\eta\}\cap A)>0$. Define $A':=\{R_{t^*}\le\nu-\eta\}\cap A$ and $\delta:=\eta^2/\beta^2$. Then $(\tilde R^\nu_{t^*})^2\ge\delta>0$ on $A'$.

Nonnegativity of $(\tilde R^\nu_u)^2\|\nabla_\theta\log\psi_u^{\theta_n}\|^2\psi_u^{\theta_n}$ gives
\[
\E\bigl[\|\hat\nabla L_{\rm score}(\theta_n)\|^2\bigr] \;\ge\; \E\!\left[\mathbf{1}_{A'}\, (\tilde R^\nu_{t^*})^2\, \|\nabla_\theta\log\psi_{t^*}^{\theta_n}\|^2\, \psi_{t^*}^{\theta_n}\right].
\]

Define
\[
V_n:=\sum_{s<t^*}\frac{\nabla_\theta\phi_s^{\theta_n}(S_s)}{1-\phi_s^{\theta_n}(S_s)},
\qquad
c:=\inf_n\inf_{\omega\in A}\|\nabla_\theta\phi_{t^*}^{\theta_n}(S_{t^*}(\omega))\|^2>0.
\]
On $A'$, condition (b) and Assumption~\ref{asm:smooth} give
\[
\prod_{s<t^*}(1-\phi_s^{\theta_n}(S_s))\ge\epsilon_0^{t^*},
\qquad
\|V_n\|^2\le t^*\sum_{s<t^*}\frac{G^2}{\epsilon_0^2}
=\frac{(t^*)^2G^2}{\epsilon_0^2}.
\]
Apply $\|u-v\|^2\ge\tfrac12\|u\|^2-\|v\|^2$ to
$\nabla_\theta\log\psi_{t^*}^{\theta_n}
=\nabla_\theta\phi_{t^*}^{\theta_n}(S_{t^*})/\phi_{t^*}^{\theta_n}(S_{t^*})-V_n$.
Using $\psi_{t^*}^{\theta_n}/\phi_{t^*}^{\theta_n}(S_{t^*})
=\prod_{s<t^*}(1-\phi_s^{\theta_n}(S_s))$, $\psi_{t^*}^{\theta_n}\le1$, and
$\E[(\tilde R_{t^*}^\nu)^2]\le\tilde M^2$ gives
\[
\E[\|\hat\nabla L_{\rm score}(\theta_n)\|^2]
\ge\frac{\delta c\epsilon_0^{t^*}}2
\E\!\left[\frac{\mathbf1_{A'}}{\phi_{t^*}^{\theta_n}(S_{t^*})}\right]
-\frac{\tilde M^2(t^*)^2G^2}{\epsilon_0^2}.
\]
Since $\phi_{t^*}^{\theta_n}(S_{t^*})\to0$ on $A'$ and $P(A')>0$,
Fatou's lemma gives
$\E[\mathbf1_{A'}/\phi_{t^*}^{\theta_n}(S_{t^*})]\to\infty$.
Therefore $\E[\|\hat\nabla L_{\rm score}(\theta_n)\|^2]\to\infty$.

For fixed $\nu$, nonnegativity of $R_t$ gives $|\tilde R_t^\nu|\le\nu^+/\beta$. Lemma~\ref{thm:variance}(a) therefore gives
\[
\E[\|\hat\nabla L(\theta)\|^2]\le (T+1)^2C_\psi^2(\nu^+)^2/\beta^2,
\]
for every $\theta\in\Theta$.

\end{proof}

\subsection*{Sigmoid policies}\label{app:sigmoid}

Under Assumptions~\ref{asm:reg},~\ref{asm:smooth}, and~\ref{asm:l2}, fix $\nu\in\R$ and a sequence $(\theta_n)\subset\Theta$. Let $B=1$.
Define $\varsigma(x):=(1+e^{-x})^{-1}$. Let $\phi_t^\theta(S_t)=\varsigma(z_t^\theta(S_t))$ with $z_t^\theta(S_t)\in\R$ for $t<T$. Assume
\[
\begin{gathered}
\sup_n\sup_{\omega\in\Omega}\max_{t<T}\|\nabla_\theta z_t^{\theta_n}(S_t(\omega))\|<\infty,\\
\phi_t^{\theta_n}(S_t)\to\phi_t^{\rm lim}\in\{0,1\}
\quad P\text{-a.s.\ for every }t<T.
\end{gathered}
\]
Then
\[
\E[\|\hat\nabla L_{\rm score}(\theta_n)\|^2]\to0,
\qquad
\E[\|\hat\nabla L(\theta_n)\|^2]\to0.
\]

\begin{proof}
Define $\tilde M:=\max_{0\le t\le T}\|\tilde R_t^\nu\|_{L^2(P)}$. Define $\tau_n:=\tau_{\theta_n}$ and $a_t^{(n)}:=\mathbf1\{U_t\le\phi_t^{\theta_n}(S_t)\}$ for $t<T$, with $U_0,\ldots,U_T$ as in~\eqref{eq:tau}.
For $t<T$, the sigmoid parameterization gives
$\nabla_\theta\log\pi_\theta(a_t\mid S_t)=(a_t-\phi_t^\theta(S_t))\nabla_\theta z_t^\theta(S_t)$.
Therefore
\[
\nabla_\theta\log\pi_\theta(\tau)
=\sum_{t=0}^{\min(\tau,T-1)}(a_t-\phi_t^\theta(S_t))\nabla_\theta z_t^\theta(S_t).
\]
Define $C_z:=\sup_n\sup_{\omega\in\Omega}\max_{t<T}\|\nabla_\theta z_t^{\theta_n}(S_t(\omega))\|<\infty$.
Since $|a_t^{(n)}-\phi_t^{\theta_n}(S_t)|\le1$, the triangle inequality gives
$\|\nabla_\theta\log\pi_{\theta_n}(\tau_n)\|\le(T+1)C_z$.

Since $(\tilde R^\nu_{\tau_n})^2\le\sum_{s=0}^T(\tilde R^\nu_s)^2$ and $\|\nabla_\theta\log\pi_{\theta_n}(\tau_n)\|\le(T+1)C_z$,
\[
\|\hat\nabla L_{\rm score}(\theta_n)\|^2 \;\le\; (1 + T)^2 C_z^2 \sum_{s=0}^T (\tilde R^\nu_s)^2 \quad \text{pointwise},
\]
with $\E[\sum_{s=0}^T (\tilde R^\nu_s)^2] \le (T+1) \tilde M^2 < \infty$, independent of $n$.

For each $t<T$, $\phi_t^{\theta_n}(S_t)\to\phi_t^{\rm lim}\in\{0,1\}$ and $U_t\ne\phi_t^{\rm lim}$ almost surely. Hence $a_t^{(n)}=\phi_t^{\rm lim}$ for all sufficiently large $n$, so $a_t^{(n)}-\phi_t^{\theta_n}(S_t)\to0$ almost surely. Since $T$ is finite and $\|\nabla_\theta z_t^{\theta_n}(S_t)\|\le C_z$, the finite score sum satisfies $\nabla_\theta\log\pi_{\theta_n}(\tau_n)\to0$ almost surely.

Hence $\|\hat\nabla L_{\rm score}(\theta_n)\|^2\to0$ almost surely. Dominated convergence with the integrable bound $(T+1)^2C_z^2\sum_{s=0}^T(\tilde R_s^\nu)^2$ gives $\E[\|\hat\nabla L_{\rm score}(\theta_n)\|^2]\to0$.

For $t<T$, $\nabla_\theta\phi_t^{\theta_n}(S_t)=\phi_t^{\theta_n}(S_t)(1-\phi_t^{\theta_n}(S_t))\nabla_\theta z_t^{\theta_n}(S_t)\to0$. The derivative formula in Remark~\ref{rem:boundary} gives $\nabla_\theta\psi_t^{\theta_n}\to0$ for $0\le t\le T$. Thus $\hat\nabla L(\theta_n)\to0$ almost surely. Since $|\tilde R_t^\nu|\le\nu^+/\beta$, the pathwise formula gives $\|\hat\nabla L(\theta_n)\|\le(T+1)C_\psi\nu^+/\beta$. Dominated convergence yields $\E[\|\hat\nabla L(\theta_n)\|^2]\to0$.

For each $t<T$, the sigmoid derivative also gives
\[
\|\nabla_\theta\phi_t^{\theta_n}(S_t)\|
\le C_z\phi_t^{\theta_n}(S_t)(1-\phi_t^{\theta_n}(S_t)).
\]
Thus $\nabla_\theta\phi_t^{\theta_n}(S_t)\to0$ wherever $\phi_t^{\theta_n}(S_t)\to0$, which excludes condition (a) of Theorem~\ref{cor:pathwise-vs-score} when the logit gradients are uniformly bounded.
\end{proof}

\subsection*{Proof of Theorem~\ref{thm:conv}}\label{app:convproof}

\begin{proof}
Let $\theta\in\Theta_{\nu_-}$ and define $q:=\nu^*(\theta)$. Let $X_1,\ldots,X_B$ be independent copies of $R_{\tau_\theta}$. Let $X_{(1)}\le\cdots\le X_{(B)}$ denote the sorted sample and define $\widehat q:=X_{(m_B)}$. Define
\[
H:=\max_{0\le t\le T}R_t,\qquad M_H:=\sqrt{T+1}\,M,\qquad
Q:=\frac{M_H}{\sqrt{1-\beta}}.
\]
Then $\E[H^2]\le M_H^2$. Assumption~\ref{asm:conv}(a) at $u=0$ gives $F_\theta(q)=\beta$. Since $0\le R_{\tau_\theta}\le H$, the quantile $q$ satisfies $(1-\beta)q^2\le\E[(R_{\tau_\theta})^2]\le M_H^2$, and hence $0\le q\le Q$. Define
\begin{equation}
\label{eq:Bzero}
B_0:=\left\lceil\max\left\{3,\frac{2}{f_-\delta_0},
1+\frac{2\beta}{f_-\delta_0}\right\}\right\rceil.
\end{equation}

Let $B\ge B_0$. For $x\in\R$, define $N_x:=\sum_{i=1}^B\mathbf1_{\{X_i\le x\}}$. The order statistic $\widehat q$ satisfies
\[
\{\widehat q>q+u\}\subseteq\{N_{q+u}\le m_B-1\},
\qquad
\{\widehat q<q-u\}\subseteq\{N_{q-u}\ge m_B\}.
\]
For $2/(f_-B)\le u\le\delta_0$, Assumption~\ref{asm:conv}(a) and the rank condition \eqref{eq:rank} give
\[
\E[N_{q+u}]-(m_B-1)\ge Bf_-u,\qquad
m_B-\E[N_{q-u}]\ge Bf_-u-1\ge Bf_-u/2.
\]
Hoeffding's inequality~\citep[Theorem 2]{hoeffding1963} for the binomial counts $N_{q+u}$ and $N_{q-u}$ therefore yields
\begin{equation}
\label{eq:localtail}
\Prob(|\widehat q-q|>u)\le2\exp\!\left(-\frac{Bf_-^2u^2}{2}\right).
\end{equation}
Integrating~\eqref{eq:localtail} over $2/(f_-B)\le u\le\delta_0$ and using $\Prob(|\widehat q-q|>u)\le1$ for $0\le u<2/(f_-B)$ gives
\begin{align}
\E\!\left[\min\{(\widehat q-q)^2,\delta_0^2\}\right]
&=\int_0^{\delta_0}2u\Prob(|\widehat q-q|>u)\,du\notag\\
&\le\frac{4}{f_-^2B^2}
+\int_0^\infty4u\exp\!\left(-\frac{Bf_-^2u^2}{2}\right)du
\le\frac{8}{f_-^2B}.
\label{eq:localmoment}
\end{align}

Define $E_+:=\{\widehat q>q+\delta_0\}$ and $E_-:=\{\widehat q<q-\delta_0\}$.
For each $i=1,\ldots,B$, define
\[
N_{q+\delta_0}^{(-i)}:=\sum_{j\ne i}\mathbf1_{\{X_j\le q+\delta_0\}}.
\]
The count $N_{q+\delta_0}^{(-i)}$ is independent of $X_i$, has distribution
$\operatorname{Bin}(B-1,F_\theta(q+\delta_0))$, and satisfies
$E_+\subseteq\{N_{q+\delta_0}^{(-i)}\le m_B-1\}$.
By \eqref{eq:Bzero}, \eqref{eq:rank}, and Assumption~\ref{asm:conv}(a),
\[
\frac{m_B-1}{B-1}\le\beta+\frac{\beta}{B-1}
\le\beta+\frac{f_-\delta_0}{2},
\qquad F_\theta(q+\delta_0)\ge\beta+f_-\delta_0.
\]
Independence of $X_i$ and $N_{q+\delta_0}^{(-i)}$ and Hoeffding's inequality give
\[
\E[X_i^2\mathbf1_{E_+}]
\le\E[X_i^2]\Prob(N_{q+\delta_0}^{(-i)}\le m_B-1)
\le M_H^2\exp\!\left(-\frac{(B-1)f_-^2\delta_0^2}{2}\right).
\]
Since $\widehat q^2\le\sum_{i=1}^B X_i^2$,
\begin{equation}
\label{eq:outerupper}
\E[\widehat q^2\mathbf1_{E_+}]
\le BM_H^2\exp\!\left(-\frac{(B-1)f_-^2\delta_0^2}{2}\right).
\end{equation}
Hoeffding's inequality applied to $N_{q+\delta_0}$ and $N_{q-\delta_0}$ gives
\begin{equation}
\label{eq:outerprob}
\Prob(E_+)\le e^{-Bf_-^2\delta_0^2/2},\qquad
\Prob(E_-)\le e^{-Bf_-^2\delta_0^2/2}.
\end{equation}
On $E_+$, $(\widehat q-q)^2\le2\widehat q^2+2Q^2$.
On $E_-$, nonnegativity of $X_1,\ldots,X_B$ gives $(\widehat q-q)^2\le Q^2$.
Combining \eqref{eq:outerupper} and \eqref{eq:outerprob} gives
\begin{equation}
\label{eq:outermoment}
\E[(\widehat q-q)^2\mathbf1_{E_+\cup E_-}]
\le(2BM_H^2+3Q^2)
\exp\!\left(-\frac{(B-1)f_-^2\delta_0^2}{2}\right).
\end{equation}
Since $B(2BM_H^2+3Q^2)e^{-(B-1)f_-^2\delta_0^2/2}$ is bounded for $B\ge B_0$,
\eqref{eq:localmoment} and \eqref{eq:outermoment} imply that there exists $C_\nu>0$,
depending only on $\beta$, $T$, $M$, $f_-$, and $\delta_0$, such that
\begin{equation}
\label{eq:mseuniform}
\E[(\widehat q-q)^2]\le C_\nu/B
\quad\text{for every }\theta\in\Theta_{\nu_-}\text{ and }B\ge B_0.
\end{equation}

Conditional on $\mathcal H_k$, the rewards $R_{\tau_{\theta_k}^{(i)}}^{(i)}$, $i=1,\ldots,B$, are independent with CDF $F_{\theta_k}$. Applying \eqref{eq:mseuniform} with $\theta=\theta_k$ therefore gives
\begin{equation}
\label{eq:conditionalquantile}
\E[(\nu_k-\nu^*(\theta_k))^2\mid\mathcal H_k]\le C_\nu/B.
\end{equation}
Taking expectation in \eqref{eq:conditionalquantile} proves \eqref{eq:quantile}.

For $\theta\in\Theta_{\nu_-}$ and $q=\nu^*(\theta)$, define
\begin{equation}
\label{eq:Wq}
W_q(\vartheta):=q-\frac1\beta\E\!\left[\sum_{t=0}^T
\psi_t^\vartheta(q-R_t)^+\right],\qquad \vartheta\in U.
\end{equation}
The representation~\eqref{eq:ru} gives $V_\beta(\vartheta)\ge W_q(\vartheta)$ for every $\vartheta\in U$ and $V_\beta(\theta)=W_q(\theta)$.
The bounds $0\le(q-R_t)^+\le q$ and $\|\nabla_\vartheta\psi_t^\vartheta\|\le C_\psi$ justify differentiating $W_q$ under the expectation.
By Assumption~\ref{asm:conv}(b), the function $V_\beta-W_q$ is differentiable on $U$ and has a minimum at $\theta\in U$. Hence
\begin{equation}
\label{eq:contact}
\nabla V_\beta(\theta)=\nabla W_q(\theta)
=-\frac1\beta\E\!\left[\sum_{t=0}^T(q-R_t)^+\nabla_\theta\psi_t^\theta\right].
\end{equation}

Define
\[
\Lambda:=\frac{(T+1)C_\psi}{\beta},\qquad
q_k:=\nu^*(\theta_k),\qquad v_k:=\nabla V_\beta(\theta_k).
\]
For each sampled trajectory $i$, define
\[
Y_k^{(i)}:=-\frac1\beta\sum_{t=0}^T(q_k-R_t^{(i)})^+
\nabla_\theta\psi_t^{\theta_k,(i)}.
\]
Then $\widehat g_k(q_k)=B^{-1}\sum_{i=1}^B Y_k^{(i)}$.
Conditional on $\mathcal H_k$, the vectors $Y_k^{(1)},\ldots,Y_k^{(B)}$ are independent.
Equation~\eqref{eq:contact} gives $\E[Y_k^{(i)}\mid\mathcal H_k]=v_k$, and
Assumption~\ref{asm:smooth} with $q_k\le Q$ gives $\|Y_k^{(i)}\|\le\Lambda Q$.
Define $\xi_k:=\widehat g_k(q_k)-v_k$. Conditional independence of $Y_k^{(1)},\ldots,Y_k^{(B)}$ gives
\begin{equation}
\label{eq:noise}
\E[\xi_k\mid\mathcal H_k]=0,\qquad
\E[\|\xi_k\|^2\mid\mathcal H_k]\le\frac{\Lambda^2Q^2}{B}.
\end{equation}
Define $\Delta_k:=\widehat g_k(\nu_k)-\widehat g_k(q_k)$.
The inequality $|(\nu_k-R_t^{(i)})^+-(q_k-R_t^{(i)})^+|\le|\nu_k-q_k|$ and
the derivative bound $\|\nabla_\theta\psi_t^{\theta_k,(i)}\|\le C_\psi$ in Assumption~\ref{asm:smooth} give $\|\Delta_k\|\le\Lambda|\nu_k-q_k|$.
Applying \eqref{eq:conditionalquantile} yields
\begin{equation}
\label{eq:thresholdmoment}
\E[\|\Delta_k\|^2\mid\mathcal H_k]\le\frac{\Lambda^2C_\nu}{B}.
\end{equation}
Define
\[
e_k:=\widehat g_k(\nu_k)-v_k,\qquad
C_g:=2\Lambda^2(Q^2+C_\nu).
\]
Since $e_k=\xi_k+\Delta_k$, the inequality $\|\xi_k+\Delta_k\|^2\le2\|\xi_k\|^2+2\|\Delta_k\|^2$,
together with \eqref{eq:noise} and \eqref{eq:thresholdmoment}, gives
\begin{equation}
\label{eq:gradienterror}
\E[\|e_k\|^2\mid\mathcal H_k]\le C_g/B.
\end{equation}

Since $\theta_k,\theta_{k+1}\in U$ and $U$ is convex, Assumption~\ref{asm:conv}(b) and $\theta_{k+1}=\theta_k+\eta(v_k+e_k)$ imply
\begin{align*}
V_\beta(\theta_{k+1})-V_\beta(\theta_k)
&\ge\eta\langle v_k,v_k+e_k\rangle
-\frac{L_V\eta^2}{2}\|v_k+e_k\|^2\\
&=\frac\eta2\bigl(\|v_k\|^2-\|e_k\|^2\bigr)
+\frac\eta2(1-L_V\eta)\|v_k+e_k\|^2.
\end{align*}
Since $1-L_V\eta\ge0$, taking conditional expectation and applying \eqref{eq:gradienterror} yields
\begin{equation}
\label{eq:ascent}
\E[V_\beta(\theta_{k+1})-V_\beta(\theta_k)\mid\mathcal H_k]
\ge\frac\eta2\|v_k\|^2-\frac{\eta C_g}{2B}.
\end{equation}
For every $\vartheta\in\Theta$, the inequalities $0\le V_\beta(\vartheta)\le\E[H]<\infty$ imply that $D_0$ is finite.
Taking expectation in \eqref{eq:ascent}, summing over $k=0,\ldots,K-1$, and applying
$\E[V_\beta(\theta_K)-V_\beta(\theta_0)]\le D_0$ gives
\[
\frac1K\sum_{k=0}^{K-1}\E[\|v_k\|^2]
\le\frac{2D_0}{\eta K}+\frac{C_g}{B},
\]
which proves \eqref{eq:stationarity}.

The quantile bound \eqref{eq:quantile} and the gradient bound \eqref{eq:stationarity} give
\begin{align*}
&\min_{0\le k<K}\E[\|\nabla V_\beta(\theta_k)\|^2+(\nu_k-\nu^*(\theta_k))^2]\\
&\qquad\le\frac1K\sum_{k=0}^{K-1}\E[\|\nabla V_\beta(\theta_k)\|^2+(\nu_k-\nu^*(\theta_k))^2]
\le\frac{2D_0}{\eta K}+\frac{C_g+C_\nu}{B}.
\end{align*}
Equation~\eqref{eq:conv} follows. Let $\epsilon>0$ and define
\begin{equation}
\label{eq:allocation}
K=\max\left\{1,\left\lceil\frac{4D_0}{\eta\epsilon}\right\rceil\right\},
\qquad
B=\max\left\{B_0,\left\lceil\frac{2(C_g+C_\nu)}{\epsilon}\right\rceil\right\}.
\end{equation}
By~\eqref{eq:allocation}, $2D_0/(\eta K)\le\epsilon/2$ and $(C_g+C_\nu)/B\le\epsilon/2$.
For fixed $\beta$ and $\eta$,
$K=O(\epsilon^{-1})$, $B=O(\epsilon^{-1})$, and $KB=O(\epsilon^{-2})$ for sufficiently small $\epsilon>0$.
\end{proof}

\section{Implementation details for the symmetric Bermudan max-call domain}
\label{app:exp-details}

We follow the LunarLander protocol of~\citet{luo2026cvarpg} with the following changes. (i) The discount $\gamma = 1$ since the max-call payoff $e^{-rt\Delta t}(\max_j S_t^j-K_{\rm strike})^+$ already discounts through the $e^{-rt\Delta t}$ factor (with $\Delta t = 1/3$ year). (ii) The feature vector is $(S_t/100,\, t/T,\, R_t)$, providing the time index and the immediate payoff suited to optimal stopping. (iii) BatchNorm is applied to the encoder. (iv) The policy logit at $t=T$ is set to $+20$. (v) Policy gradients are clipped to a maximum norm of $1.0$. (vi) The IQL batch size of MIX scales linearly with horizon $T$ from the LunarLander value in~\citet{luo2024mixture}: $2 \times 10^5$ transitions per update at $T = 500$ becomes $3{,}600$ at $T = 9$. The IQL replay buffer holds $10^6$ transitions. Tables~\ref{tab:domain},~\ref{tab:protocol}, and~\ref{tab:method-hp} specify the domain components, training protocol, and per-method hyperparameters, respectively.

\paragraph{Policy and quantile implementation.}
The encoder has two Linear--BatchNorm--ReLU blocks of width 128, followed by a $128\to64\to1$ policy head with ReLU and sigmoid output. During training, BatchNorm computes statistics across the $B(T+1)$ sampled state-time pairs. Dropout is disabled. Rewards are divided by the standard deviation of the $B(T+1)$ payoffs in the first batch, and the normalization scale is held fixed throughout training.
The implementation estimates the reward quantile using the order statistic of rank $\max\{1,\lfloor\beta_kB\rfloor\}$ from the sampled stopping rewards. Floating-point schedule arithmetic gives hold-phase ranks $3,1,1$ at target levels $0.2,0.1,0.05$ with $B=20$. Two annealing updates at target $0.2$ use ranks $7$ and $5$ instead of $8$ and $6$. With annealing disabled, the separate branch uses rank $4$ at $0.2$. The loss is the stop-gradient surrogate in Algorithm~\ref{alg:dios}, with $\log\psi_t$ computed from log-sigmoid sums and $\nu_k$ held fixed during differentiation.

\paragraph{Compute.}
Experiments ran on a server with eight NVIDIA RTX A6000 GPUs, with each job assigned to one GPU. Some jobs shared a GPU. Training times include trajectory simulation, gradient updates, logging, and six periodic evaluations of $5{,}000$ paths. Model creation and the final $20{,}000$-path evaluation are excluded. RET-CAP times also exclude training and evaluation of the policies used to set the minimum cap. The RET-CAP CVaR values come from reruns with recomputed minimum caps and retuned hyperparameters, while the tables retain the training times measured with the earlier minimum caps and hyperparameters.

\paragraph{Evaluation.} Each final policy is evaluated on $20{,}000$ paths with seed $10000+s$ for training seed $s\in\{0,\ldots,9\}$. Methods share evaluation paths for each seed. Evaluation uses BatchNorm running statistics and stops at the first $t$ with $\phi_t^\theta(S_t)>0.5$. At $\beta=0.2$, CVaR is the mean of exactly the lowest $4{,}000$ rewards. Means and population standard deviations are over ten training seeds.

\begin{table}[ht]
\centering
\caption{Bermudan max-call domain parameters.}
\label{tab:domain}
\small
\setlength{\tabcolsep}{6pt}
\begin{tabular}{ll}
\toprule
Parameter & Value \\
\midrule
Strike $K_{\rm strike}$ & $100$ \\
Risk-free rate $r$ & $0.05$ \\
Dividend yield $\delta$ & $0.10$ \\
Volatility $\sigma$ & $0.20$ (shared across all $d$ underlyings) \\
Asset correlation $\rho$ & $0$ (independent) \\
Maturity & $3$ years \\
Last time index $T$ & $9$ ($10$ candidate times including $t=0$) \\
Time step $\Delta t$ & $1/3$ year \\
Payoff & $g_t(S_t) = e^{-rt\Delta t}(\max_j S_t^j-K_{\rm strike})^+$ \\
Policy input features & $(S_t/100,\, t/T,\, R_t)$ \\
$d$-grid & $\{20, 50, 100, 200\}$ \\
$s_0$-grid & $\{90, 100, 110\}$ \\
Configurations & $12$ \\
\bottomrule
\end{tabular}
\end{table}

\begin{table}[ht]
\centering
\caption{Shared training protocol.}
\label{tab:protocol}
\small
\setlength{\tabcolsep}{6pt}
\begin{tabular}{ll}
\toprule
Item & Value \\
\midrule
Trajectories per update $B$ & $20$ \\
Outer iterations $K$ & $6{,}000$ \\
Risk level $\beta$ & $0.2$ \\
Discount $\gamma$ & $1$ \\
Seeds & $10$ \\
Optimizer & Adam \\
Policy LR sweep set & $\{7\mathrm{e}{-4},\, 5\mathrm{e}{-4},\, 3\mathrm{e}{-4},\, 1\mathrm{e}{-4},\, 7\mathrm{e}{-5}\}$ \\
Value LR multiplier sweep set & $\{1,2,5,10\}$ times the policy learning rate \\
Advantage normalization sweep set & $\{$yes, no$\}$ where applicable \\
Hidden width & $128$ \\
Hidden layers & $2$ \\
Activation & ReLU \\
Encoder normalization & BatchNorm \\
Terminal logit overwrite & $+20$ at $t=T$ \\
Multi-step $\lambda$ (CVaR-VaR, RET-CAP only) & $0.95$ \\
Gradient clip (max-norm) & $1.0$ \\
Final evaluation batch & $20{,}000$ \\
\bottomrule
\end{tabular}
\end{table}

\begin{table}[ht]
\centering
\caption{Per-method best hyperparameters from the learning-rate sweep on the representative configuration $d = 100$, $s_0 = 100$, applied uniformly to all twelve configurations.}
\label{tab:method-hp}
\footnotesize
\setlength{\tabcolsep}{3pt}
\begin{tabular}{>{\raggedright\arraybackslash}p{1.6cm}>{\raggedright\arraybackslash}p{1.4cm}>{\raggedright\arraybackslash}p{0.9cm}>{\raggedright\arraybackslash}p{1.0cm}>{\raggedright\arraybackslash}p{7.0cm}}
\toprule
Method & Policy LR & Critic LR mult. & Norm. adv. & Method-specific \\
\midrule
DIOS (ours) & $7\mathrm{e}{-5}$ & --- & --- & 3-phase $\beta_k$-schedule with warmup fraction $\rho_w = 0.25$, anneal fraction $\rho_a = 0.50$, hold fraction $0.25$: $\beta_k = 1$ for the first $\rho_w K$ iterations, linear $1 \to 0.2$ over the next $\rho_a K$, $\beta_k = 0.2$ thereafter \\
CVaR-VaR & $1\mathrm{e}{-4}$ & $1\times$ & no & $N_{\mathrm{q}} = 10$ atoms in the quantile critic; tradeoff weight $\omega = 0.5$ throughout training, selected from the three $\omega$-schedules of~\citet{luo2026cvarpg} (Maze: $\omega \equiv 0.5$; LunarLander: $\omega = 0.5$ for first 40\% then linear decay to 0; Pendulum: $\omega = 0.5$ for first 25\% then 0) \\
CVaR-PG & $3\mathrm{e}{-4}$ & --- & no & trajectory-level score-function estimator with no value critic~\citep{tamar2015cvar} \\
RET-CAP (cap-anchor) & $1\mathrm{e}{-4}$ & $5\times$ & yes & minimum cap $q^*$ computed from DIOS evaluation quantiles as in Appendix~\ref{app:retcap-anchor}; advantage-shaping weight $0.8$ \\
RET-CAP (RN-anchor) & $7\mathrm{e}{-5}$ & $10\times$ & yes & risk-neutral anchor variant of Appendix~\ref{app:retcap-rn}; advantage-shaping weight $0.8$ \\
MIX & $7\mathrm{e}{-4}$ & --- & --- & IQL component: LR $1\mathrm{e}{-4}$, update interval $50$, expectile $0.8$, temperature $\beta_{\mathrm{IQL}} = 1$, batch size $3{,}600$ per update, replay buffer capacity $10^6$ transitions \\
\bottomrule
\end{tabular}
\end{table}

\paragraph{Paired comparisons.}\label{app:paired}
Table~\ref{tab:paired-beta02} reports two-sided paired $t$ and exact Wilcoxon tests over training seeds, using shared evaluation paths per seed and no multiple-comparison correction. DIOS has a higher mean CVaR than CVaR-VaR in all 24 max-call configurations, with paired $t$-test $p<0.05$ against CVaR-VaR in 2 symmetric and 3 asymmetric configurations.

\begin{table}[ht]
\centering
\caption{Seed-paired difference DIOS $-$ CVaR-VaR in CVaR$_{0.2}$ (10 seeds, shared evaluation paths per seed). $t$: two-sided paired $t$; W: two-sided exact Wilcoxon signed-rank; no multiple-comparison correction.}
\label{tab:paired-beta02}
\small
\setlength{\tabcolsep}{5pt}
\begin{tabular}{cc ccccc}
\toprule
$d$ & $s_0$ & $\bar\Delta$ & 95\% CI & $t$ $p$ & W $p$ & $p_t<0.05$ \\
\midrule
\multicolumn{7}{l}{Symmetric} \\
20 & 90 & $+0.103$ & [-0.064, +0.271] & 0.197 & 0.625 & no \\
20 & 100 & $+0.113$ & [-0.086, +0.311] & 0.23 & 0.492 & no \\
20 & 110 & $+0.073$ & [-0.204, +0.350] & 0.565 & 1 & no \\
50 & 90 & $+0.007$ & [-0.081, +0.095] & 0.857 & 1 & no \\
50 & 100 & $+0.088$ & [-0.010, +0.186] & 0.0736 & 0.232 & no \\
50 & 110 & $+0.151$ & [-0.078, +0.381] & 0.169 & 0.232 & no \\
100 & 90 & $+0.142$ & [+0.010, +0.274] & 0.0374 & 0.0273 & yes \\
100 & 100 & $+0.032$ & [-0.065, +0.129] & 0.472 & 0.322 & no \\
100 & 110 & $+0.125$ & [+0.017, +0.232] & 0.0274 & 0.00977 & yes \\
200 & 90 & $+0.033$ & [-0.071, +0.136] & 0.494 & 0.232 & no \\
200 & 100 & $+0.025$ & [-0.051, +0.102] & 0.467 & 0.375 & no \\
200 & 110 & $+0.003$ & [-0.090, +0.096] & 0.94 & 0.77 & no \\
\midrule
\multicolumn{7}{l}{Asymmetric} \\
20 & 90 & $+0.812$ & --- & 0.00523 & 0.00586 & yes \\
20 & 100 & $+0.427$ & --- & 0.195 & 0.193 & no \\
20 & 110 & $+0.857$ & --- & 0.0187 & 0.0273 & yes \\
50 & 90 & $+0.451$ & --- & 0.0522 & 0.0371 & no \\
50 & 100 & $+0.568$ & --- & 0.0828 & 0.0488 & no \\
50 & 110 & $+1.418$ & --- & 0.018 & 0.00977 & yes \\
100 & 90 & $+0.745$ & --- & 0.114 & 0.193 & no \\
100 & 100 & $+0.058$ & --- & 0.817 & 0.922 & no \\
100 & 110 & $+0.654$ & --- & 0.281 & 0.922 & no \\
200 & 90 & $+0.265$ & --- & 0.184 & 0.232 & no \\
200 & 100 & $+0.526$ & --- & 0.172 & 0.432 & no \\
200 & 110 & $+0.071$ & --- & 0.844 & 0.77 & no \\
\bottomrule
\end{tabular}
\end{table}

\section{Additional baseline experiments on the symmetric Bermudan max-call domain}
\label{app:additional}

\subsection{PCVaR-PG}
\label{app:pcvar}

Table~\ref{tab:pcvar_main} reports PCVaR-PG~\citep{kim2024pcvar} on all twelve symmetric configurations. At $s_0=110$, all ten seeds stop immediately and receive reward $10$. At $s_0\in\{90,100\}$, the mean stopping time ranges from $1.0$ to $3.4$.

\begin{table}[t]
\centering
\caption{PCVaR-PG~\citep{kim2024pcvar} on the Bermudan max-call domain across all 12 configurations, with selected policy learning rate $3\mathrm{e}{-4}$ from a sweep on $(d, s_0)=(100, 100)$. CVaR$_{0.2}$ mean $\pm$ std over 10 seeds, per-run wall time, and mean stopping time $\bar\tau \in [0, T]$ with $T = 9$.}
\label{tab:pcvar_main}
\small
\begin{tabular}{ccccc}
\toprule
$d$ & $s_0$ & CVaR$_{0.2}$ & wall (s) & $\bar\tau$ \\
\midrule
\multirow{3}{*}{20}
& 90  & $2.741  \pm 0.638$ & $29.5$ & $2.4$ \\
& 100 & $11.407 \pm 0.254$ & $30.1$ & $3.4$ \\
& 110 & $10.000 \pm 0.000$ & $27.5$ & $0.0$ \\
\midrule
\multirow{3}{*}{50}
& 90  & $6.461  \pm 0.149$ & $33.6$ & $3.2$ \\
& 100 & $18.024 \pm 0.050$ & $32.5$ & $3.1$ \\
& 110 & $10.000 \pm 0.000$ & $26.4$ & $0.0$ \\
\midrule
\multirow{3}{*}{100}
& 90  & $8.265  \pm 4.137$ & $31.6$ & $1.7$ \\
& 100 & $22.267 \pm 0.053$ & $32.9$ & $2.1$ \\
& 110 & $10.000 \pm 0.000$ & $30.9$ & $0.0$ \\
\midrule
\multirow{3}{*}{200}
& 90  & $12.401 \pm 4.134$ & $31.5$ & $1.3$ \\
& 100 & $26.196 \pm 0.042$ & $26.5$ & $1.0$ \\
& 110 & $10.000 \pm 0.000$ & $30.1$ & $0.0$ \\
\bottomrule
\end{tabular}
\end{table}

\subsection{RET-CAP minimum cap and the choice of anchor}
\label{app:retcap-anchor}

RET-CAP updates a return cap during training and constrains the cap to be at least an externally supplied value $q^*$~\citep{mead2025retcap}. In the max-call experiments, we compute the empirical $\beta$-quantile of the evaluation rewards for each trained DIOS policy. We average the quantiles over ten seeds and round to three decimal places to obtain $q^*$. The tables label the DIOS-derived minimum cap ``cap-anchor''. Appendix~\ref{app:retcap-rn} compares the DIOS-derived minimum cap with a minimum cap computed from separately trained risk-neutral policies. The choice of an evaluation VaR as the minimum cap is part of our experimental protocol.

\subsection{RET-CAP with risk-neutral anchor}
\label{app:retcap-rn}

We additionally set the RET-CAP minimum cap from the evaluation quantiles of separately trained risk-neutral (RN) policies on the symmetric domain. Table~\ref{tab:retcap_rn} compares the DIOS-derived and risk-neutral minimum caps, with hyperparameters selected separately for each choice. The asymmetric counterpart is reported in Appendix~\ref{app:asym-retcap-rn}.

\begin{table}[t]
\centering
\caption{RET-CAP with two anchor choices on the Bermudan max-call domain across all 12 configurations: the default cap-anchor (used in our main results) and the risk-neutral anchor variant. CVaR$_{0.2}$ mean $\pm$ std over 10 seeds; $\Delta$ is RN-anchor CVaR minus cap-anchor CVaR; $\bar\tau_{\text{RN}}$ is the mean stopping time under the RN-anchor variant ($T = 9$).}
\label{tab:retcap_rn}
\small
\begin{tabular}{cccccc}
\toprule
$d$ & $s_0$ & cap CVaR$_{0.2}$ & RN CVaR$_{0.2}$ & $\Delta$ & $\bar\tau_{\text{RN}}$ \\
\midrule
\multirow{3}{*}{20} &  90 & $18.003\,\pm\,0.290$ & $16.447\,\pm\,0.395$ & $-1.556$ & $3.2$ \\
 & 100 & $29.677\,\pm\,0.295$ & $27.107\,\pm\,0.377$ & $-2.570$ & $2.8$ \\
 & 110 & $41.616\,\pm\,0.366$ & $37.921\,\pm\,0.459$ & $-3.696$ & $2.5$ \\
\midrule
\multirow{3}{*}{50} &  90 & $32.810\,\pm\,0.225$ & $30.483\,\pm\,0.378$ & $-2.327$ & $3.9$ \\
 & 100 & $45.913\,\pm\,0.278$ & $42.976\,\pm\,0.705$ & $-2.936$ & $3.7$ \\
 & 110 & $59.115\,\pm\,0.379$ & $55.346\,\pm\,0.580$ & $-3.768$ & $3.4$ \\
\midrule
\multirow{3}{*}{100} &  90 & $44.323\,\pm\,0.359$ & $42.148\,\pm\,0.655$ & $-2.175$ & $4.6$ \\
 & 100 & $58.956\,\pm\,0.390$ & $56.053\,\pm\,0.712$ & $-2.904$ & $4.4$ \\
 & 110 & $73.084\,\pm\,0.600$ & $69.170\,\pm\,0.910$ & $-3.914$ & $4.1$ \\
\midrule
\multirow{3}{*}{200} &  90 & $56.105\,\pm\,0.423$ & $53.672\,\pm\,0.402$ & $-2.433$ & $5.0$ \\
 & 100 & $71.396\,\pm\,0.410$ & $68.896\,\pm\,0.620$ & $-2.500$ & $4.9$ \\
 & 110 & $86.622\,\pm\,0.795$ & $84.240\,\pm\,0.589$ & $-2.381$ & $4.8$ \\
\bottomrule
\end{tabular}
\end{table}

\subsection{DIOS ablations}
\label{app:ablations}

We compare variants of two design choices in DIOS: the soft-penalty CVaR objective and the $\beta$-anneal schedule.

\paragraph{Risk-neutral baseline.} We train the same DIOS architecture under the risk-neutral expectation $J_{\rm RN}(\theta) = \E[\sum_t \psi_t^\theta R_t]$, replacing the soft-penalty objective with expected reward and omitting quantile estimation. We evaluate the risk-neutral policy and DIOS using CVaR$_{0.2}$ on the same evaluation paths for each seed. Table~\ref{tab:rn_baseline_sym} reports CVaR$_{0.2}$ for the risk-neutral policy and DIOS, together with the difference, risk-neutral minus DIOS.

\begin{table}[t]
\centering
\caption{Risk-neutral DIOS baseline vs DIOS on the symmetric domain. CVaR$_{0.2}$ at evaluation, mean over 10 seeds; $\Delta = $ RN $-$ DIOS. RN training drops the soft penalty and trains the same architecture under risk-neutral expectation; the resulting policy is then evaluated against the CVaR$_{0.2}$ target.}
\label{tab:rn_baseline_sym}
\small
\setlength{\tabcolsep}{6pt}
\begin{tabular}{cc ccc}
\toprule
$d$ & $s_0$ & DIOS & RN & $\Delta$ \\
\midrule
\multirow{3}{*}{20} &  90 & $18.17$ & $9.59$ & $-8.58$ \\
 & 100 & $29.96$ & $19.92$ & $-10.05$ \\
 & 110 & $41.97$ & $30.56$ & $-11.41$ \\
\midrule
\multirow{3}{*}{50} &  90 & $33.22$ & $25.79$ & $-7.43$ \\
 & 100 & $46.80$ & $38.25$ & $-8.56$ \\
 & 110 & $60.39$ & $50.70$ & $-9.69$ \\
\midrule
\multirow{3}{*}{100} &  90 & $45.00$ & $38.35$ & $-6.65$ \\
 & 100 & $59.84$ & $52.19$ & $-7.65$ \\
 & 110 & $74.67$ & $66.03$ & $-8.64$ \\
\midrule
\multirow{3}{*}{200} &  90 & $57.14$ & $50.88$ & $-6.26$ \\
 & 100 & $73.26$ & $66.12$ & $-7.14$ \\
 & 110 & $89.38$ & $81.36$ & $-8.02$ \\
\bottomrule
\end{tabular}
\end{table}

\paragraph{$\beta$-anneal ablation.} We disable both the warmup and anneal phases, setting $\rho_w = \rho_a = 0$ so that $\beta_k = \beta_{\rm tgt} = 0.2$ from iteration $k = 0$. Table~\ref{tab:frac0_ablation_sym} reports the comparison on the symmetric domain. Without the schedule, all ten seeds collapse at $(d,s_0)=(20,110)$, and two of ten seeds collapse in each configuration with $d\in\{100,200\}$. The full-schedule and no-schedule runs use hold-phase quantile ranks three and four, respectively.

\begin{table}[t]
\centering
\caption{$\beta$-anneal ablation on the symmetric domain: DIOS with full schedule ($\rho_w = 0.25$, $\rho_a = 0.50$) vs no schedule ($\rho_w = \rho_a = 0$, $\beta = \beta_{\rm tgt} = 0.2$ from iteration $0$). CVaR$_{0.2}$ mean over 10 seeds; coll reports the number of collapsed seeds out of ten, with collapse defined by $\bar\tau<0.5$.}
\label{tab:frac0_ablation_sym}
\small
\setlength{\tabcolsep}{6pt}
\begin{tabular}{cc ccc}
\toprule
$d$ & $s_0$ & DIOS & No schedule & coll \\
\midrule
\multirow{3}{*}{20} &  90 & $18.17$ & $17.65$ & $0/10$ \\
 & 100 & $29.96$ & $29.57$ & $0/10$ \\
 & 110 & $41.97$ & $10.00$ & $10/10$ \\
\midrule
\multirow{3}{*}{50} &  90 & $33.22$ & $33.01$ & $0/10$ \\
 & 100 & $46.80$ & $46.51$ & $0/10$ \\
 & 110 & $60.39$ & $60.02$ & $0/10$ \\
\midrule
\multirow{3}{*}{100} &  90 & $45.00$ & $35.88$ & $2/10$ \\
 & 100 & $59.84$ & $47.68$ & $2/10$ \\
 & 110 & $74.67$ & $61.51$ & $2/10$ \\
\midrule
\multirow{3}{*}{200} &  90 & $57.14$ & $45.66$ & $2/10$ \\
 & 100 & $73.26$ & $58.55$ & $2/10$ \\
 & 110 & $89.38$ & $73.47$ & $2/10$ \\
\bottomrule
\end{tabular}
\end{table}

\subsection{CVaR-VaR critic ablation}
\label{app:cvarvar-critic}

CVaR-VaR~\citep{luo2026cvarpg} uses a state-dependent quantile critic to track quantile levels and compute multi-step advantages for the VaR-gradient term. Our implementation uses a separate $128$-unit MLP trained by multi-step quantile regression at each iteration, adding a critic forward pass and optimizer step to each policy update.

The no-critic variant replaces the critic prediction at every state with the preceding batch's mean stopping reward after normalization. At the first update, the baseline is the mean unnormalized terminal payoff in the batch. The no-critic variant also holds the tradeoff weight $\omega$ at $0.5$ for the first 40\% of training and then decreases $\omega$ linearly to zero. The main CVaR-VaR baseline uses $\omega=0.5$ throughout.

In every configuration in Table~\ref{tab:cvarvar-critic-sym}, the no-critic variant has shorter recorded training time and lower mean CVaR than the main CVaR-VaR baseline.

\begin{table}[ht]
\centering
\caption{CVaR-VaR critic ablation on the symmetric Bermudan max-call domain. ``Full'' is our CVaR-VaR implementation with the state-dependent quantile critic; ``no-critic'' replaces the critic with a global mean baseline, with the changed tradeoff schedule described in the text. CVaR$_{0.2}$ mean $\pm$ std and wall time in seconds, both over 10 seeds. $\Delta_\text{CVaR}$ = no-critic mean minus full mean. Speedup = wall(full) / wall(no-critic).}
\label{tab:cvarvar-critic-sym}
\small
\setlength{\tabcolsep}{4pt}
\begin{tabular}{cc|cc|cc|cc}
\toprule
$d$ & $s_0$ & full CVaR & wall & no-critic CVaR & wall & $\Delta_\text{CVaR}$ & speedup \\
\midrule
\multirow{3}{*}{20}
  & 90 & $18.068 \pm 0.385$ & $54.7$ & $11.675 \pm 0.682$ & $29.1$ & $-6.393$ & $1.88\times$ \\

  & 100 & $29.850 \pm 0.441$ & $51.0$ & $22.475 \pm 0.508$ & $28.4$ & $-7.375$ & $1.80\times$ \\

  & 110 & $41.893 \pm 0.587$ & $51.6$ & $33.078 \pm 0.561$ & $27.5$ & $-8.815$ & $1.88\times$ \\
\midrule
\multirow{3}{*}{50}
  & 90 & $33.213 \pm 0.175$ & $48.8$ & $27.879 \pm 0.449$ & $26.7$ & $-5.334$ & $1.83\times$ \\

  & 100 & $46.715 \pm 0.187$ & $58.2$ & $40.830 \pm 0.863$ & $29.0$ & $-5.885$ & $2.00\times$ \\

  & 110 & $60.241 \pm 0.311$ & $59.6$ & $53.875 \pm 0.690$ & $27.5$ & $-6.366$ & $2.17\times$ \\
\midrule
\multirow{3}{*}{100}
  & 90 & $44.859 \pm 0.195$ & $52.6$ & $40.908 \pm 0.698$ & $28.5$ & $-3.951$ & $1.85\times$ \\

  & 100 & $59.805 \pm 0.148$ & $49.9$ & $54.537 \pm 0.630$ & $27.5$ & $-5.268$ & $1.81\times$ \\

  & 110 & $74.543 \pm 0.202$ & $59.5$ & $68.702 \pm 0.699$ & $29.7$ & $-5.841$ & $2.01\times$ \\
\midrule
\multirow{3}{*}{200}
  & 90 & $57.110 \pm 0.117$ & $59.0$ & $52.955 \pm 0.385$ & $29.3$ & $-4.155$ & $2.01\times$ \\

  & 100 & $73.233 \pm 0.126$ & $58.6$ & $68.589 \pm 0.493$ & $29.7$ & $-4.644$ & $1.97\times$ \\

  & 110 & $89.375 \pm 0.099$ & $57.0$ & $84.030 \pm 0.460$ & $29.0$ & $-5.345$ & $1.97\times$ \\
\bottomrule
\end{tabular}
\end{table}

\section{Experiments on the asymmetric Bermudan max-call domain}
\label{app:asym-experiments}

\subsection{Setting and main results}
\label{app:asym}

We further evaluate the methods of \S\ref{sec:exp} on the heterogeneous-volatility variant of the Bermudan max-call, which sets $\sigma_j=0.10+j/(2d)$ for $j=1,\ldots,d$ in place of the shared $\sigma = 0.20$. The remaining domain parameters and the shared training protocol are given in Tables~\ref{tab:domain} and~\ref{tab:protocol}. Per-method best learning rates from a re-run sweep on $(d, s_0) = (100, 100)$, applied uniformly across all twelve configurations, are listed in Table~\ref{tab:method-hp-asym}. In the symmetric and asymmetric domains, DIOS uses $\rho_w=0.25$, $\rho_a=0.50$, and $\eta=7\mathrm{e}{-5}$, leaving the final $25\%$ of training at the target tail level. DIOS attains the highest mean CVaR and shortest recorded training time in all twelve configurations.

\begin{table}[ht]
\centering
\caption{Per-method best learning rates on the asymmetric Bermudan max-call domain, from a re-run sweep on the representative configuration $d = 100$, $s_0 = 100$, applied uniformly to all twelve configurations without re-tuning. Method-specific items not shown are inherited from Table~\ref{tab:method-hp}.}
\label{tab:method-hp-asym}
\small
\setlength{\tabcolsep}{6pt}
\begin{tabular}{>{\raggedright\arraybackslash}p{2.2cm}>{\raggedright\arraybackslash}p{1.6cm}>{\raggedright\arraybackslash}p{1.2cm}>{\raggedright\arraybackslash}p{1.6cm}>{\raggedright\arraybackslash}p{5.2cm}}
\toprule
Method & Policy LR & Critic LR mult. & Norm. adv. & Notes \\
\midrule
DIOS (ours) & $7\mathrm{e}{-5}$ & --- & --- & same as symmetric (Table~\ref{tab:method-hp}) \\
CVaR-VaR & $3\mathrm{e}{-4}$ & $2\times$ & no & Maze schedule ($\omega \equiv 0.5$), as in Table~\ref{tab:method-hp} \\
CVaR-PG & $3\mathrm{e}{-4}$ & --- & no & selected by our learning-rate sweep \\
MIX & $7\mathrm{e}{-4}$ & --- & --- & same as symmetric (Table~\ref{tab:method-hp}) \\
RET-CAP (cap-anchor) & $1\mathrm{e}{-4}$ & $5\times$ & yes & DIOS-derived anchor as in \S\ref{sec:exp} \\
RET-CAP (RN-anchor) & $7\mathrm{e}{-4}$ & $1\times$ & yes & risk-neutral anchor variant of Appendix~\ref{app:asym-retcap-rn} \\
\bottomrule
\end{tabular}
\end{table}

\begin{table}[ht]
\centering
\caption{Asymmetric Bermudan max-call ($\sigma_j=0.10+j/(2d)$). CVaR$_{0.2}\,\uparrow$ on top, wall (s) with collapse rate $\downarrow$ on bottom; bold marks the per-row best. The RET-CAP timing caveat in Table~\ref{tab:main} also applies here.}
\label{tab:asym}
\footnotesize
\setlength{\tabcolsep}{1.0pt}
\begin{tabular}{cc l ccccc}
\toprule
$d$ & $s_0$ & & DIOS (ours) & CVaR-VaR & CVaR-PG & MIX & RET-CAP \\
\midrule
 \multirow{6}{*}{$20$} & \multirow{2}{*}{$90$}
 & CVaR\,$\uparrow$ & $\mathbf{49.572}\,\pm\,0.587$ & $48.761\,\pm\,0.831$ & $43.825\,\pm\,14.656$ & $38.857\,\pm\,19.481$ & $49.028\,\pm\,0.572$ \\
 & & wall\,(coll)\,$\downarrow$ & $\mathbf{12.4}$\,(0/10) & $55.8$\,(0/10) & $15.7$\,(1/10) & $21.7$\,(2/10) & $25.6$\,(0/10) \\
\cmidrule(l){2-8}
  & \multirow{2}{*}{$100$}
 & CVaR\,$\uparrow$ & $\mathbf{64.962}\,\pm\,0.753$ & $64.535\,\pm\,1.296$ & $62.877\,\pm\,2.676$ & $58.096\,\pm\,19.378$ & $64.277\,\pm\,0.782$ \\
 & & wall\,(coll)\,$\downarrow$ & $\mathbf{11.8}$\,(0/10) & $62.1$\,(0/10) & $15.6$\,(0/10) & $21.7$\,(1/10) & $24.8$\,(0/10) \\
\cmidrule(l){2-8}
  & \multirow{2}{*}{$110$}
 & CVaR\,$\uparrow$ & $\mathbf{80.425}\,\pm\,0.822$ & $79.568\,\pm\,1.027$ & $58.401\,\pm\,31.707$ & $50.893\,\pm\,33.456$ & $79.806\,\pm\,0.696$ \\
 & & wall\,(coll)\,$\downarrow$ & $\mathbf{12.1}$\,(0/10) & $52.2$\,(0/10) & $15.9$\,(3/10) & $23.1$\,(4/10) & $24.3$\,(0/10) \\
\midrule
 \multirow{6}{*}{$50$} & \multirow{2}{*}{$90$}
 & CVaR\,$\uparrow$ & $\mathbf{95.068}\,\pm\,0.532$ & $94.617\,\pm\,0.907$ & $94.485\,\pm\,0.771$ & $94.037\,\pm\,1.703$ & $93.887\,\pm\,1.039$ \\
 & & wall\,(coll)\,$\downarrow$ & $\mathbf{11.7}$\,(0/10) & $63.5$\,(0/10) & $15.2$\,(0/10) & $21.7$\,(0/10) & $25.3$\,(0/10) \\
\cmidrule(l){2-8}
  & \multirow{2}{*}{$100$}
 & CVaR\,$\uparrow$ & $\mathbf{115.542}\,\pm\,0.587$ & $114.974\,\pm\,1.131$ & $114.553\,\pm\,0.969$ & $103.134\,\pm\,34.408$ & $114.633\,\pm\,0.905$ \\
 & & wall\,(coll)\,$\downarrow$ & $\mathbf{11.8}$\,(0/10) & $61.5$\,(0/10) & $15.8$\,(0/10) & $23.0$\,(1/10) & $24.2$\,(0/10) \\
\cmidrule(l){2-8}
  & \multirow{2}{*}{$110$}
 & CVaR\,$\uparrow$ & $\mathbf{135.981}\,\pm\,0.629$ & $134.563\,\pm\,1.741$ & $132.866\,\pm\,3.655$ & $122.375\,\pm\,37.468$ & $134.756\,\pm\,0.827$ \\
 & & wall\,(coll)\,$\downarrow$ & $\mathbf{11.8}$\,(0/10) & $54.7$\,(0/10) & $14.8$\,(0/10) & $22.3$\,(1/10) & $26.9$\,(0/10) \\
\midrule
 \multirow{6}{*}{$100$} & \multirow{2}{*}{$90$}
 & CVaR\,$\uparrow$ & $\mathbf{139.929}\,\pm\,0.828$ & $139.184\,\pm\,1.232$ & $125.096\,\pm\,41.709$ & $125.357\,\pm\,41.805$ & $139.242\,\pm\,0.935$ \\
 & & wall\,(coll)\,$\downarrow$ & $\mathbf{11.8}$\,(0/10) & $53.5$\,(0/10) & $15.9$\,(1/10) & $22.6$\,(1/10) & $26.7$\,(0/10) \\
\cmidrule(l){2-8}
  & \multirow{2}{*}{$100$}
 & CVaR\,$\uparrow$ & $\mathbf{165.289}\,\pm\,0.862$ & $165.230\,\pm\,1.085$ & $148.721\,\pm\,49.581$ & $148.515\,\pm\,49.525$ & $164.846\,\pm\,1.192$ \\
 & & wall\,(coll)\,$\downarrow$ & $\mathbf{12.2}$\,(0/10) & $59.8$\,(0/10) & $14.2$\,(1/10) & $23.7$\,(1/10) & $26.6$\,(0/10) \\
\cmidrule(l){2-8}
  & \multirow{2}{*}{$110$}
 & CVaR\,$\uparrow$ & $\mathbf{190.629}\,\pm\,1.000$ & $189.975\,\pm\,2.070$ & $171.654\,\pm\,53.899$ & $172.614\,\pm\,54.209$ & $189.273\,\pm\,1.732$ \\
 & & wall\,(coll)\,$\downarrow$ & $\mathbf{12.9}$\,(0/10) & $55.4$\,(0/10) & $15.2$\,(1/10) & $23.7$\,(1/10) & $24.4$\,(0/10) \\
\midrule
 \multirow{6}{*}{$200$} & \multirow{2}{*}{$90$}
 & CVaR\,$\uparrow$ & $\mathbf{195.053}\,\pm\,1.137$ & $194.788\,\pm\,1.104$ & $171.552\,\pm\,58.019$ & $175.488\,\pm\,58.503$ & $193.957\,\pm\,1.301$ \\
 & & wall\,(coll)\,$\downarrow$ & $\mathbf{13.4}$\,(0/10) & $62.3$\,(0/10) & $15.7$\,(1/10) & $22.9$\,(1/10) & $25.9$\,(0/10) \\
\cmidrule(l){2-8}
  & \multirow{2}{*}{$100$}
 & CVaR\,$\uparrow$ & $\mathbf{226.624}\,\pm\,1.351$ & $226.099\,\pm\,1.399$ & $198.697\,\pm\,67.013$ & $203.191\,\pm\,67.779$ & $224.970\,\pm\,1.756$ \\
 & & wall\,(coll)\,$\downarrow$ & $\mathbf{13.6}$\,(0/10) & $54.2$\,(0/10) & $14.3$\,(1/10) & $24.5$\,(1/10) & $24.9$\,(0/10) \\
\cmidrule(l){2-8}
  & \multirow{2}{*}{$110$}
 & CVaR\,$\uparrow$ & $\mathbf{258.131}\,\pm\,1.444$ & $258.060\,\pm\,1.147$ & $228.388\,\pm\,73.747$ & $232.825\,\pm\,74.287$ & $255.473\,\pm\,2.795$ \\
 & & wall\,(coll)\,$\downarrow$ & $\mathbf{12.5}$\,(0/10) & $51.0$\,(0/10) & $15.8$\,(1/10) & $23.2$\,(1/10) & $27.1$\,(0/10) \\
\bottomrule
\end{tabular}
\end{table}

\subsection{PCVaR-PG}
\label{app:asym-pcvar}

Table~\ref{tab:pcvar_main_asym} reports PCVaR-PG on the asymmetric domain. At $s_0=110$, every seed stops at $t=0$ and receives reward $10$. At $s_0\in\{90,100\}$, mean stopping times are at most $4.5$, compared with $T=9$. The standard deviation of CVaR$_{0.2}$ reaches $14.185$ at $(d,s_0)=(200,90)$.

\begin{table}[t]
\centering
\caption{PCVaR-PG~\citep{kim2024pcvar} on the asymmetric Bermudan max-call domain across all 12 configurations, with selected policy learning rate $3\mathrm{e}{-4}$ from a sweep on $(d, s_0)=(100, 100)$. CVaR$_{0.2}$ mean $\pm$ std over 10 seeds, per-run wall time, and mean stopping time $\bar\tau \in [0, T]$ with $T = 9$.}
\label{tab:pcvar_main_asym}
\small
\begin{tabular}{ccccc}
\toprule
$d$ & $s_0$ & CVaR$_{0.2}$ & wall (s) & $\bar\tau$ \\
\midrule
\multirow{3}{*}{20}
&  90 & $7.205  \pm 4.774$ & $29.8$ & $1.4$ \\
& 100 & $16.519 \pm 8.289$ & $33.1$ & $2.4$ \\
& 110 & $10.000 \pm 0.000$ & $29.8$ & $0.0$ \\
\midrule
\multirow{3}{*}{50}
&  90 & $23.860 \pm 0.295$ & $34.4$ & $3.1$ \\
& 100 & $37.740 \pm 1.456$ & $33.4$ & $4.5$ \\
& 110 & $10.000 \pm 0.000$ & $30.2$ & $0.0$ \\
\midrule
\multirow{3}{*}{100}
&  90 & $28.266 \pm 14.134$ & $33.6$ & $2.5$ \\
& 100 & $50.220 \pm 0.219$  & $31.9$ & $2.5$ \\
& 110 & $10.000 \pm 0.000$  & $28.7$ & $0.0$ \\
\midrule
\multirow{3}{*}{200}
&  90 & $42.551 \pm 14.185$ & $30.3$ & $0.9$ \\
& 100 & $50.752 \pm 25.377$ & $28.8$ & $0.9$ \\
& 110 & $10.000 \pm 0.000$  & $29.1$ & $0.0$ \\
\bottomrule
\end{tabular}
\end{table}

\subsection{RET-CAP with risk-neutral anchor}
\label{app:asym-retcap-rn}

Table~\ref{tab:retcap_rn_asym} compares RET-CAP with DIOS-derived and risk-neutral minimum caps on the asymmetric domain. Hyperparameters are selected separately for each minimum cap.

\begin{table}[t]
\centering
\caption{RET-CAP with two anchor choices on the asymmetric Bermudan max-call domain across all 12 configurations: cap-anchor (default, used in main results, Table~\ref{tab:asym}) and risk-neutral anchor variant. CVaR$_{0.2}$ mean $\pm$ std over 10 seeds; $\Delta$ is RN-anchor CVaR minus cap-anchor CVaR; $\bar\tau_{\text{RN}}$ is the mean stopping time under the RN-anchor variant ($T = 9$).}
\label{tab:retcap_rn_asym}
\small
\begin{tabular}{cccccc}
\toprule
$d$ & $s_0$ & cap CVaR$_{0.2}$ & RN CVaR$_{0.2}$ & $\Delta$ & $\bar\tau_{\text{RN}}$ \\
\midrule
\multirow{3}{*}{20} &  90 & $49.028\,\pm\,0.572$ & $43.606\,\pm\,2.596$ & $-5.422$ & $2.7$ \\
 & 100 & $64.277\,\pm\,0.782$ & $55.169\,\pm\,4.124$ & $-9.109$ & $2.4$ \\
 & 110 & $79.806\,\pm\,0.696$ & $66.089\,\pm\,3.564$ & $-13.717$ & $2.1$ \\
\midrule
\multirow{3}{*}{50} &  90 & $93.887\,\pm\,1.039$ & $84.134\,\pm\,2.521$ & $-9.753$ & $3.4$ \\
 & 100 & $114.633\,\pm\,0.905$ & $102.160\,\pm\,2.142$ & $-12.473$ & $3.3$ \\
 & 110 & $134.756\,\pm\,0.827$ & $119.518\,\pm\,1.263$ & $-15.237$ & $3.1$ \\
\midrule
\multirow{3}{*}{100} &  90 & $139.242\,\pm\,0.935$ & $125.880\,\pm\,2.753$ & $-13.362$ & $3.9$ \\
 & 100 & $164.846\,\pm\,1.192$ & $151.648\,\pm\,3.372$ & $-13.198$ & $4.0$ \\
 & 110 & $189.273\,\pm\,1.732$ & $173.946\,\pm\,2.998$ & $-15.328$ & $3.9$ \\
\midrule
\multirow{3}{*}{200} &  90 & $193.957\,\pm\,1.301$ & $179.846\,\pm\,6.094$ & $-14.111$ & $4.6$ \\
 & 100 & $224.970\,\pm\,1.756$ & $212.469\,\pm\,3.740$ & $-12.500$ & $4.6$ \\
 & 110 & $255.473\,\pm\,2.795$ & $239.484\,\pm\,5.234$ & $-15.989$ & $4.5$ \\
\bottomrule
\end{tabular}
\end{table}

\subsection{DIOS ablations}
\label{app:asym-ablations}

We repeat the risk-neutral-objective and annealing comparisons from Appendix~\ref{app:ablations} on the asymmetric domain.

\paragraph{Risk-neutral baseline.} Table~\ref{tab:rn_baseline_asym} reports mean CVaR$_{0.2}$ for the risk-neutral policy and DIOS, together with the difference, risk-neutral minus DIOS.

\begin{table}[t]
\centering
\caption{Risk-neutral DIOS baseline vs DIOS on the asymmetric domain. CVaR$_{0.2}$ at evaluation, mean over 10 seeds; $\Delta = $ RN $-$ DIOS.}
\label{tab:rn_baseline_asym}
\small
\setlength{\tabcolsep}{6pt}
\begin{tabular}{cc ccc}
\toprule
$d$ & $s_0$ & DIOS & RN & $\Delta$ \\
\midrule
\multirow{3}{*}{20} &  90 & $49.57$ & $24.53$ & $-25.05$ \\
 & 100 & $64.96$ & $36.63$ & $-28.33$ \\
 & 110 & $80.42$ & $48.86$ & $-31.56$ \\
\midrule
\multirow{3}{*}{50} &  90 & $95.07$ & $66.16$ & $-28.91$ \\
 & 100 & $115.54$ & $83.08$ & $-32.46$ \\
 & 110 & $135.98$ & $100.00$ & $-35.98$ \\
\midrule
\multirow{3}{*}{100} &  90 & $139.93$ & $108.71$ & $-31.22$ \\
 & 100 & $165.29$ & $130.37$ & $-34.92$ \\
 & 110 & $190.63$ & $152.01$ & $-38.62$ \\
\midrule
\multirow{3}{*}{200} &  90 & $195.05$ & $160.66$ & $-34.40$ \\
 & 100 & $226.62$ & $188.10$ & $-38.52$ \\
 & 110 & $258.13$ & $215.50$ & $-42.63$ \\
\bottomrule
\end{tabular}
\end{table}

\paragraph{$\beta$-anneal ablation.} Table~\ref{tab:frac0_ablation_asym} reports the comparison on the asymmetric domain. The no-schedule variant has more collapsed seeds, defined by $\bar\tau<0.5$, than DIOS with the full schedule.

\begin{table}[t]
\centering
\caption{$\beta$-anneal ablation on the asymmetric domain: DIOS with full schedule ($\rho_w = 0.25$, $\rho_a = 0.50$) vs no schedule ($\rho_w = \rho_a = 0$, $\beta = \beta_{\rm tgt} = 0.2$ from iteration $0$). CVaR$_{0.2}$ mean over 10 seeds; coll reports the number of collapsed seeds out of ten, with collapse defined by $\bar\tau<0.5$.}
\label{tab:frac0_ablation_asym}
\small
\setlength{\tabcolsep}{6pt}
\begin{tabular}{cc ccc}
\toprule
$d$ & $s_0$ & DIOS & No schedule & coll \\
\midrule
\multirow{3}{*}{20} &  90 & $49.57$ & $49.27$ & $0/10$ \\
 & 100 & $64.96$ & $64.53$ & $0/10$ \\
 & 110 & $80.42$ & $10.00$ & $10/10$ \\
\midrule
\multirow{3}{*}{50} &  90 & $95.07$ & $94.20$ & $0/10$ \\
 & 100 & $115.54$ & $114.58$ & $0/10$ \\
 & 110 & $135.98$ & $134.88$ & $0/10$ \\
\midrule
\multirow{3}{*}{100} &  90 & $139.93$ & $111.36$ & $2/10$ \\
 & 100 & $165.29$ & $131.55$ & $2/10$ \\
 & 110 & $190.63$ & $153.71$ & $2/10$ \\
\midrule
\multirow{3}{*}{200} &  90 & $195.05$ & $154.34$ & $2/10$ \\
 & 100 & $226.62$ & $179.28$ & $2/10$ \\
 & 110 & $258.13$ & $206.25$ & $2/10$ \\
\bottomrule
\end{tabular}
\end{table}

\subsection{CVaR-VaR critic ablation}
\label{app:cvarvar-critic-asym}

We repeat the CVaR-VaR critic ablation of Appendix~\ref{app:cvarvar-critic} on the asymmetric domain. The no-critic variant uses the batch-mean baseline defined in Appendix~\ref{app:cvarvar-critic} and decreases $\omega$ to zero after holding $\omega=0.5$ for the first 40\% of training. Table~\ref{tab:cvarvar-critic-asym} compares the no-critic variant with CVaR-VaR using a critic and $\omega=0.5$ throughout.

\begin{table}[ht]
\centering
\caption{CVaR-VaR critic ablation on the asymmetric Bermudan max-call domain. ``Full'' is our CVaR-VaR implementation with the state-dependent quantile critic; ``no-critic'' replaces the critic with a global mean baseline, with the changed tradeoff schedule described in the text. CVaR$_{0.2}$ mean $\pm$ std and wall time in seconds, both over 10 seeds. $\Delta_\text{CVaR}$ = no-critic mean minus full mean. Speedup = wall(full) / wall(no-critic).}
\label{tab:cvarvar-critic-asym}
\small
\setlength{\tabcolsep}{4pt}
\begin{tabular}{cc|cc|cc|cc}
\toprule
$d$ & $s_0$ & full CVaR & wall & no-critic CVaR & wall & $\Delta_\text{CVaR}$ & speedup \\
\midrule
\multirow{3}{*}{20}
  & 90 & $48.761 \pm 0.831$ & $55.8$ & $29.581 \pm 0.949$ & $28.2$ & $-19.180$ & $1.98\times$ \\

  & 100 & $64.535 \pm 1.296$ & $62.1$ & $43.363 \pm 3.072$ & $27.7$ & $-21.172$ & $2.24\times$ \\

  & 110 & $79.568 \pm 1.027$ & $52.2$ & $55.406 \pm 1.149$ & $26.2$ & $-24.162$ & $1.99\times$ \\
\midrule
\multirow{3}{*}{50}
  & 90 & $94.617 \pm 0.907$ & $63.5$ & $71.847 \pm 1.509$ & $28.1$ & $-22.770$ & $2.26\times$ \\

  & 100 & $114.974 \pm 1.131$ & $61.5$ & $91.802 \pm 1.437$ & $29.1$ & $-23.172$ & $2.11\times$ \\

  & 110 & $134.563 \pm 1.741$ & $54.7$ & $106.933 \pm 0.933$ & $27.2$ & $-27.630$ & $2.01\times$ \\
\midrule
\multirow{3}{*}{100}
  & 90 & $139.184 \pm 1.232$ & $53.5$ & $116.141 \pm 2.076$ & $28.1$ & $-23.043$ & $1.91\times$ \\

  & 100 & $165.230 \pm 1.085$ & $59.8$ & $140.633 \pm 3.651$ & $29.4$ & $-24.597$ & $2.03\times$ \\

  & 110 & $189.975 \pm 2.070$ & $55.4$ & $162.248 \pm 2.673$ & $28.7$ & $-27.727$ & $1.93\times$ \\
\midrule
\multirow{3}{*}{200}
  & 90 & $194.788 \pm 1.104$ & $62.3$ & $170.293 \pm 2.300$ & $29.0$ & $-24.495$ & $2.15\times$ \\

  & 100 & $226.099 \pm 1.399$ & $54.2$ & $197.651 \pm 2.369$ & $26.6$ & $-28.448$ & $2.03\times$ \\

  & 110 & $258.060 \pm 1.147$ & $51.0$ & $228.504 \pm 2.271$ & $29.4$ & $-29.556$ & $1.73\times$ \\
\bottomrule
\end{tabular}
\end{table}

\clearpage
\section{Smaller tail levels and quantile smoothing}
\label{app:smallbeta}
Table~\ref{tab:small-beta} reports the symmetric max-call grid at $\beta=0.1$ and $0.05$. Raw DIOS uses the current batch order statistic as the quantile estimate. For DIOS with an exponential moving average (EMA), denoted $\nu$-EMA, let $\widehat\nu_k$ be the batch order statistic and define
\[
\nu_0=\widehat\nu_0,\qquad \nu_k=(1-\alpha)\nu_{k-1}+\alpha\widehat\nu_k\quad(k\ge1).
\]
The smoothing is applied from the beginning of training, including warmup and annealing. The weight $\alpha$ multiplies the new batch statistic. We select $\alpha=0.9$ at $\beta=0.1$ on $(d,s_0)=(200,100)$ and $\alpha=0.1$ at $\beta=0.05$ on $(d,s_0)=(20,90)$, using seeds 0--2. The final comparison uses ten seeds. DIOS retains learning rate $7\times10^{-5}$. CVaR-VaR uses learning rate $10^{-4}$, value multiplier one, and fixed tradeoff weight $\omega=0.5$. CVaR-PG uses learning rate $3\times10^{-4}$.

DIOS with EMA exceeds CVaR-VaR and CVaR-PG in mean CVaR in all 12 configurations at each tail level. Two-sided paired $t$-tests over ten seeds give $p<0.05$ against CVaR-VaR in 6/12 configurations at $\beta=0.1$ and 12/12 at $\beta=0.05$. The $p$-values are not adjusted for multiple comparisons. Raw DIOS has a higher mean CVaR than CVaR-VaR and CVaR-PG in 1 of 12 configurations at $\beta=0.1$ and 8 of 12 configurations at $\beta=0.05$. The improvement from smoothing is observed with batch size $B=20$.

\begin{table}[ht]
\centering
\caption{Lower tail levels, symmetric max-call. CVaR$_\beta$ mean $\pm$ SD over 10 seeds. DIOS raw: batch order statistic; DIOS $\nu$-EMA: $\nu_k = (1-\alpha)\nu_{k-1} + \alpha\hat\nu_k$. Bold: highest mean.}
\label{tab:small-beta}
\small
\setlength{\tabcolsep}{5pt}
\begin{tabular}{cc cccc}
\toprule
$d$ & $s_0$ & DIOS (raw) & DIOS ($\nu$-EMA) & CVaR-VaR & CVaR-PG \\
\midrule
\multicolumn{6}{l}{$\beta = 0.1$ ($\nu$-EMA $\alpha = 0.9$)} \\
20 & 90 & $14.117\,\pm\,0.489$ & $\mathbf{14.759}\,\pm\,0.282$ & $14.548\,\pm\,0.392$ & $3.974\,\pm\,4.598$ \\
20 & 100 & $25.467\,\pm\,0.576$ & $\mathbf{26.202}\,\pm\,0.359$ & $25.823\,\pm\,0.573$ & $10.850\,\pm\,11.412$ \\
20 & 110 & $37.086\,\pm\,0.857$ & $\mathbf{37.854}\,\pm\,0.393$ & $37.547\,\pm\,0.571$ & $12.678\,\pm\,8.033$ \\
50 & 90 & $29.108\,\pm\,0.462$ & $\mathbf{29.410}\,\pm\,0.214$ & $29.266\,\pm\,0.374$ & $13.527\,\pm\,11.258$ \\
50 & 100 & $42.299\,\pm\,0.492$ & $\mathbf{42.676}\,\pm\,0.206$ & $42.334\,\pm\,0.532$ & $24.711\,\pm\,15.589$ \\
50 & 110 & $55.428\,\pm\,0.559$ & $\mathbf{55.943}\,\pm\,0.233$ & $55.531\,\pm\,0.456$ & $36.452\,\pm\,16.694$ \\
100 & 90 & $40.863\,\pm\,0.204$ & $\mathbf{41.125}\,\pm\,0.171$ & $40.922\,\pm\,0.374$ & $14.152\,\pm\,16.300$ \\
100 & 100 & $55.181\,\pm\,0.229$ & $\mathbf{55.567}\,\pm\,0.212$ & $55.226\,\pm\,0.480$ & $21.000\,\pm\,22.708$ \\
100 & 110 & $69.559\,\pm\,0.282$ & $\mathbf{70.014}\,\pm\,0.234$ & $69.813\,\pm\,0.306$ & $32.822\,\pm\,24.598$ \\
200 & 90 & $52.868\,\pm\,0.161$ & $\mathbf{53.127}\,\pm\,0.136$ & $52.813\,\pm\,0.435$ & $27.531\,\pm\,23.773$ \\
200 & 100 & $68.620\,\pm\,0.124$ & $\mathbf{68.808}\,\pm\,0.132$ & $68.634\,\pm\,0.289$ & $36.662\,\pm\,31.060$ \\
200 & 110 & $84.236\,\pm\,0.158$ & $\mathbf{84.506}\,\pm\,0.162$ & $84.395\,\pm\,0.225$ & $49.657\,\pm\,33.571$ \\
\midrule
\multicolumn{6}{l}{$\beta = 0.05$ ($\nu$-EMA $\alpha = 0.1$)} \\
20 & 90 & $7.472\,\pm\,0.892$ & $\mathbf{12.123}\,\pm\,0.246$ & $9.363\,\pm\,2.778$ & $1.536\,\pm\,2.131$ \\
20 & 100 & $17.404\,\pm\,1.203$ & $\mathbf{23.239}\,\pm\,0.295$ & $21.396\,\pm\,1.225$ & $5.057\,\pm\,6.297$ \\
20 & 110 & $27.787\,\pm\,1.374$ & $\mathbf{34.698}\,\pm\,0.302$ & $29.226\,\pm\,2.444$ & $15.150\,\pm\,6.519$ \\
50 & 90 & $22.412\,\pm\,0.707$ & $\mathbf{26.456}\,\pm\,0.244$ & $24.768\,\pm\,1.108$ & $9.611\,\pm\,7.502$ \\
50 & 100 & $34.889\,\pm\,1.015$ & $\mathbf{39.333}\,\pm\,0.285$ & $33.346\,\pm\,1.786$ & $19.106\,\pm\,10.825$ \\
50 & 110 & $47.367\,\pm\,1.224$ & $\mathbf{52.280}\,\pm\,0.316$ & $44.453\,\pm\,2.490$ & $30.110\,\pm\,11.272$ \\
100 & 90 & $34.293\,\pm\,0.509$ & $\mathbf{38.086}\,\pm\,0.175$ & $31.031\,\pm\,1.267$ & $10.986\,\pm\,12.276$ \\
100 & 100 & $47.926\,\pm\,0.461$ & $\mathbf{52.198}\,\pm\,0.203$ & $41.660\,\pm\,2.179$ & $17.358\,\pm\,18.225$ \\
100 & 110 & $61.674\,\pm\,0.515$ & $\mathbf{66.321}\,\pm\,0.217$ & $53.797\,\pm\,1.910$ & $28.809\,\pm\,19.661$ \\
200 & 90 & $46.731\,\pm\,0.635$ & $\mathbf{49.903}\,\pm\,0.190$ & $37.470\,\pm\,2.399$ & $16.877\,\pm\,17.550$ \\
200 & 100 & $61.650\,\pm\,0.554$ & $\mathbf{65.238}\,\pm\,0.213$ & $50.787\,\pm\,2.648$ & $23.700\,\pm\,24.174$ \\
200 & 110 & $76.879\,\pm\,0.710$ & $\mathbf{80.592}\,\pm\,0.239$ & $61.297\,\pm\,2.082$ & $35.638\,\pm\,26.123$ \\
\bottomrule
\end{tabular}
\end{table}

\FloatBarrier

\section{House selling}
\label{app:house}
Offers $X_1,\ldots,X_9$ arrive independently with mean 50. The seller observes one offer at each decision point $j=1,\ldots,9$, can accept only the current offer, and must accept by $j=9$. Accepting at $j$ gives $X_j-cj$, so the first observed offer already incurs cost $c$. There is no stopping decision before the first offer. Relabeling the offer times as $t=0,\ldots,8$ gives the time convention of Section~\ref{sec:setup}. The policy observes the normalized current offer, time, and immediate payoff.

The lognormal grid uses $\log X_j\sim N(\log50-\sigma^2/2,\sigma^2)$ with $\sigma\in\{0.3,0.5,0.7,0.9\}$. The Pareto grid uses shape $a\in\{6,4,3,2.5\}$ and scale $50(a-1)/a$. Each family is crossed with $c\in\{1,2,3\}$, giving 24 cells. Finite second moments hold in both families. Negative rewards are possible in the lognormal grid, as discussed in Appendix~\ref{app:exp-verif}.

For every method, we select the policy learning rate from the sweep set in Appendix~\ref{app:exp-details}. Selection uses $(\sigma,c)=(0.5,2)$ for lognormal offers and $(a,c)=(3,2)$ for Pareto offers. Selection and final experiments use seeds 0--9. Other settings follow the max-call protocol. RET-CAP uses the seed-averaged DIOS evaluation VaR, rounded to six decimals, as the minimum cap, and the selected RET-CAP policy learning rates are $10^{-4}$ for lognormal offers and $7\times10^{-5}$ for Pareto offers. Tables~\ref{tab:house-selling} and~\ref{tab:house-pareto} report the complete results.

RET-CAP has the largest mean only for $(\sigma,c)=(0.9,2)$, exceeding DIOS by approximately $0.004$, with paired $t$-test $p=0.99$. An alternative anchor from a policy trained with $\beta_k\equiv1$ gives RET-CAP $41.155\pm0.649$ at $(\sigma,c)=(0.9,2)$, compared with DIOS $41.729\pm0.974$. The policy used to compute the house-selling alternative minimum cap is trained at $\beta_k=1$ with rewards truncated at the batch maximum of sampled stopping rewards. The max-call risk-neutral policies instead maximize expected payoff directly.

\begin{table}[ht]
\centering
\caption{House selling: lognormal offers (9 offers, $\beta = 0.2$). CVaR$_{0.2}$ mean $\pm$ SD over 10 seeds (population SD); wall time (s) below. Bold: highest mean. RET-CAP uses the DIOS-derived anchor.}
\label{tab:house-selling}
\footnotesize
\setlength{\tabcolsep}{2.5pt}
\begin{tabular}{l l ccccc}
\toprule
cell & & DIOS & CVaR-VaR & CVaR-PG & RET-CAP & MIX \\
\midrule
ln0.3\_c1 & CVaR & $\mathbf{49.877 \pm 0.225}$ & $49.651\,\pm\,0.437$ & $49.503\,\pm\,0.523$ & $48.653\,\pm\,0.328$ & $49.029\,\pm\,0.756$ \\
 & wall & 12.3 & 49.3 & 15.4 & 23.5 & 23.3 \\
ln0.3\_c2 & CVaR & $\mathbf{45.492 \pm 0.306}$ & $45.148\,\pm\,0.382$ & $45.123\,\pm\,0.464$ & $44.843\,\pm\,0.394$ & $45.067\,\pm\,0.544$ \\
 & wall & 12.4 & 48.6 & 15.1 & 23.2 & 23.4 \\
ln0.3\_c3 & CVaR & $\mathbf{41.661 \pm 0.395}$ & $40.771\,\pm\,0.590$ & $40.829\,\pm\,0.600$ & $41.392\,\pm\,0.249$ & $40.507\,\pm\,0.651$ \\
 & wall & 12.0 & 55.1 & 15.2 & 23.0 & 23.3 \\
ln0.5\_c1 & CVaR & $\mathbf{50.647 \pm 0.453}$ & $50.048\,\pm\,0.943$ & $49.940\,\pm\,0.645$ & $49.959\,\pm\,0.412$ & $49.797\,\pm\,0.685$ \\
 & wall & 12.4 & 56.0 & 15.4 & 23.8 & 23.6 \\
ln0.5\_c2 & CVaR & $\mathbf{46.141 \pm 0.253}$ & $45.446\,\pm\,0.526$ & $45.351\,\pm\,0.611$ & $45.649\,\pm\,0.268$ & $45.396\,\pm\,0.562$ \\
 & wall & 11.9 & 57.9 & 15.5 & 23.6 & 23.9 \\
ln0.5\_c3 & CVaR & $\mathbf{41.528 \pm 0.472}$ & $40.439\,\pm\,0.827$ & $40.616\,\pm\,0.819$ & $41.206\,\pm\,0.502$ & $40.432\,\pm\,0.728$ \\
 & wall & 12.8 & 57.6 & 15.6 & 23.4 & 24.6 \\
ln0.7\_c1 & CVaR & $\mathbf{49.790 \pm 0.387}$ & $48.318\,\pm\,1.920$ & $48.444\,\pm\,1.578$ & $49.220\,\pm\,0.447$ & $49.317\,\pm\,0.620$ \\
 & wall & 12.3 & 56.8 & 15.3 & 23.6 & 24.2 \\
ln0.7\_c2 & CVaR & $\mathbf{44.933 \pm 0.307}$ & $44.103\,\pm\,0.962$ & $43.890\,\pm\,1.153$ & $44.604\,\pm\,0.449$ & $43.768\,\pm\,1.343$ \\
 & wall & 12.8 & 57.7 & 15.8 & 23.6 & 24.3 \\
ln0.7\_c3 & CVaR & $\mathbf{40.137 \pm 0.461}$ & $38.671\,\pm\,1.283$ & $38.799\,\pm\,1.184$ & $40.068\,\pm\,0.600$ & $38.696\,\pm\,1.124$ \\
 & wall & 11.9 & 53.3 & 15.0 & 24.2 & 23.6 \\
ln0.9\_c1 & CVaR & $\mathbf{47.129 \pm 0.454}$ & $45.557\,\pm\,1.703$ & $45.994\,\pm\,2.246$ & $46.764\,\pm\,0.620$ & $46.524\,\pm\,1.272$ \\
 & wall & 12.3 & 55.4 & 15.4 & 23.4 & 23.9 \\
ln0.9\_c2 & CVaR & $41.729\,\pm\,0.974$ & $40.881\,\pm\,0.974$ & $41.152\,\pm\,1.269$ & $\mathbf{41.732 \pm 0.580}$ & $39.683\,\pm\,1.740$ \\
 & wall & 12.5 & 54.2 & 15.4 & 23.8 & 23.9 \\
ln0.9\_c3 & CVaR & $\mathbf{37.309 \pm 0.269}$ & $35.500\,\pm\,2.221$ & $35.377\,\pm\,1.260$ & $37.029\,\pm\,0.498$ & $35.201\,\pm\,1.923$ \\
 & wall & 12.0 & 53.7 & 15.2 & 23.4 & 23.1 \\
\bottomrule
\end{tabular}
\end{table}


\begin{table}[ht]
\centering
\caption{House selling: Pareto offers (9 offers, $\beta = 0.2$). CVaR$_{0.2}$ mean $\pm$ SD over 10 seeds (population SD); wall time (s) below. Bold: highest mean. RET-CAP uses the DIOS-derived anchor.}
\label{tab:house-pareto}
\footnotesize
\setlength{\tabcolsep}{2.5pt}
\begin{tabular}{l l ccccc}
\toprule
cell & & DIOS & CVaR-VaR & CVaR-PG & RET-CAP & MIX \\
\midrule
pa6\_c1 & CVaR & $\mathbf{45.780 \pm 0.122}$ & $45.514\,\pm\,0.520$ & $45.322\,\pm\,0.303$ & $43.934\,\pm\,0.205$ & $45.387\,\pm\,0.366$ \\
 & wall & 12.0 & 58.3 & 15.3 & 22.2 & 23.2 \\
pa6\_c2 & CVaR & $\mathbf{42.202 \pm 0.140}$ & $41.958\,\pm\,0.550$ & $41.318\,\pm\,0.608$ & $40.674\,\pm\,0.367$ & $41.772\,\pm\,0.342$ \\
 & wall & 12.4 & 54.7 & 15.0 & 22.8 & 23.3 \\
pa6\_c3 & CVaR & $\mathbf{39.699 \pm 0.269}$ & $39.414\,\pm\,0.016$ & $38.714\,\pm\,0.477$ & $39.418\,\pm\,0.010$ & $39.123\,\pm\,0.625$ \\
 & wall & 12.0 & 63.8 & 14.9 & 21.2 & 21.4 \\
pa4\_c1 & CVaR & $\mathbf{45.279 \pm 0.171}$ & $45.088\,\pm\,0.286$ & $44.738\,\pm\,0.696$ & $43.688\,\pm\,0.300$ & $44.737\,\pm\,0.452$ \\
 & wall & 12.1 & 57.2 & 15.2 & 22.7 & 22.9 \\
pa4\_c2 & CVaR & $\mathbf{41.014 \pm 0.189}$ & $40.959\,\pm\,0.383$ & $40.590\,\pm\,0.394$ & $40.185\,\pm\,0.175$ & $40.454\,\pm\,0.379$ \\
 & wall & 12.1 & 58.0 & 15.3 & 23.0 & 23.1 \\
pa4\_c3 & CVaR & $\mathbf{37.748 \pm 0.241}$ & $37.081\,\pm\,0.699$ & $36.984\,\pm\,0.476$ & $37.405\,\pm\,0.197$ & $37.036\,\pm\,0.704$ \\
 & wall & 12.6 & 56.1 & 14.9 & 22.7 & 23.1 \\
pa3\_c1 & CVaR & $\mathbf{44.265 \pm 0.268}$ & $43.830\,\pm\,0.734$ & $43.909\,\pm\,0.322$ & $42.520\,\pm\,0.443$ & $44.092\,\pm\,0.294$ \\
 & wall & 12.1 & 56.6 & 15.5 & 23.0 & 23.2 \\
pa3\_c2 & CVaR & $\mathbf{39.655 \pm 0.200}$ & $39.423\,\pm\,0.225$ & $39.449\,\pm\,0.232$ & $38.974\,\pm\,0.346$ & $39.254\,\pm\,0.419$ \\
 & wall & 12.5 & 54.0 & 15.4 & 22.8 & 23.2 \\
pa3\_c3 & CVaR & $\mathbf{35.861 \pm 0.277}$ & $34.902\,\pm\,0.826$ & $35.036\,\pm\,0.734$ & $35.482\,\pm\,0.208$ & $35.229\,\pm\,0.517$ \\
 & wall & 12.5 & 55.5 & 15.1 & 22.9 & 23.1 \\
pa2.5\_c1 & CVaR & $\mathbf{42.991 \pm 0.206}$ & $42.556\,\pm\,0.466$ & $42.536\,\pm\,0.686$ & $41.532\,\pm\,0.438$ & $42.863\,\pm\,0.372$ \\
 & wall & 12.3 & 51.8 & 15.5 & 23.4 & 23.4 \\
pa2.5\_c2 & CVaR & $\mathbf{38.166 \pm 0.325}$ & $38.112\,\pm\,0.268$ & $37.812\,\pm\,0.329$ & $37.816\,\pm\,0.407$ & $36.430\,\pm\,2.975$ \\
 & wall & 12.8 & 51.5 & 15.9 & 24.1 & 24.6 \\
pa2.5\_c3 & CVaR & $\mathbf{34.007 \pm 0.476}$ & $33.188\,\pm\,0.928$ & $33.484\,\pm\,0.977$ & $33.831\,\pm\,0.276$ & $33.161\,\pm\,1.193$ \\
 & wall & 12.4 & 51.4 & 15.4 & 24.5 & 23.6 \\
\bottomrule
\end{tabular}
\end{table}

\FloatBarrier

\section{Network width}
\label{app:width}
Table~\ref{tab:width} compares CVaR, parameter counts, and training times for encoder widths 32, 64, 128, and 256 on symmetric $(d,s_0)=(100,100)$ at $\beta=0.2$. Mean CVaR ranges from $59.749$ to $59.838$. For the dense feedforward policy used here, a forward/backward pass costs $O(B(T+1)P)$ for $P$ parameters. Sorting the sampled rewards costs $O(B\log B)$. DIOS adds no critic network.

\begin{table}[ht]
\centering
\caption{Network width, symmetric max-call $(d, s_0) = (100, 100)$, $\beta = 0.2$. CVaR$_{0.2}$ mean $\pm$ SD over 10 seeds.}
\label{tab:width}
\small
\setlength{\tabcolsep}{5pt}
\begin{tabular}{cccc}
\toprule
hidden width & parameters & CVaR$_{0.2}$ & wall (s) \\
\midrule
32 & 5,025 & $59.785\,\pm\,0.190$ & 13.5 \\
64 & 13,121 & $59.812\,\pm\,0.201$ & 13.0 \\
128 & 38,529 & $59.838\,\pm\,0.160$ & 12.7 \\
256 & 126,209 & $59.749\,\pm\,0.110$ & 13.1 \\
\bottomrule
\end{tabular}
\end{table}

\FloatBarrier

\clearpage
\section{An optimized certainty equivalent example}
\label{app:oce}
For a nondecreasing concave utility $u$, define the optimized certainty equivalent~\citep{bental2007} by $\operatorname{OCE}_u(X):=\sup_\nu\{\nu+\E[u(X-\nu)]\}$. If $u(R_t-\nu)$ is integrable for every $t$ and every fixed $\nu$, exchanging the two suprema as in Proposition~\ref{thm:envelope} gives
\[
\sup_{\tau\in\cT}\operatorname{OCE}_u(R_\tau)
=\sup_{\nu\in\R}\left\{\nu+\sup_{\tau\in\cT}\E[u(R_\tau-\nu)]\right\}.
\]
For fixed $\nu$, the inner optimization is an expected-reward stopping problem with reward $u(R_t-\nu)$. The outer optimization selects $\nu$. For each fixed $\nu$, the Snell envelope gives an optimizer of the inner problem. If the outer supremum is attained at $\nu^*$, a stopping time that optimizes the inner problem at $\nu^*$ also optimizes the OCE stopping problem.
Define the truncated quadratic utility by
\[
u(x)=\begin{cases}x-x^2/2,&x\le1,\\1/2,&x>1,\end{cases}
\qquad u'(x)=(1-x)^+.
\]
For $X\in L^2$, any maximizing $\nu$ satisfies $\E[(1-X+\nu)^+]=1$. If $X\le1+\E[X]$ almost surely, then $\operatorname{OCE}_u(X)=\E[X]-\tfrac12\Var(X)$. The risk measure $-\operatorname{OCE}_u$ is convex but generally noncoherent. For nonnegative $R_t$, the transformed reward satisfies $u(-\nu)\le u(R_t-\nu)\le1/2$ at each fixed $\nu$, so $u(R_t-\nu)$ gives bounded weights in the pathwise gradient.

\clearpage
\section*{NeurIPS Paper Checklist}
\begin{enumerate}

\item {\bf Claims}
    \item[] Question: Do the main claims made in the abstract and introduction accurately reflect the paper's contributions and scope?
    \item[] Answer: \answerYes{}
    \item[] Justification: The abstract and introduction distinguish structural results, conditional stationarity guarantees, and empirical comparisons.


\item {\bf Limitations}
    \item[] Question: Does the paper discuss the limitations of the work performed by the authors?
    \item[] Answer: \answerYes{}
    \item[] Justification: Section~\ref{sec:conclusion} discusses limitations. Appendix~\ref{app:exp-verif} checks payoff moments and explains the policy-derivative and convergence assumptions.


\item {\bf Theory assumptions and proofs}
    \item[] Question: For each theoretical result, does the paper provide the full set of assumptions and a complete (and correct) proof?
    \item[] Answer: \answerYes{}
    \item[] Justification: Sections~\ref{sec:dios} and~\ref{sec:convergence} state the assumptions. Appendix~\ref{app:proofs} contains the proofs.


\item {\bf Experimental result reproducibility}
    \item[] Question: Does the paper fully disclose all the information needed to reproduce the main experimental results of the paper to the extent that it affects the main claims and/or conclusions of the paper (regardless of whether the code and data are provided or not)?
    \item[] Answer: \answerYes{}
    \item[] Justification: The domains, policy architecture, threshold ranks, optimizer settings, evaluation protocol, and per-method tuning are described in the appendices.


\item {\bf Open access to data and code}
    \item[] Question: Does the paper provide open access to the data and code, with sufficient instructions to faithfully reproduce the main experimental results, as described in supplemental material?
    \item[] Answer: \answerNo{}
    \item[] Justification: The manuscript does not provide a public code repository. Implementation and evaluation details are given in Appendix~\ref{app:exp-details}.


\item {\bf Experimental setting/details}
    \item[] Question: Does the paper specify all the training and test details (e.g., data splits, hyperparameters, how they were chosen, type of optimizer) necessary to understand the results?
    \item[] Answer: \answerYes{}
    \item[] Justification: Appendix~\ref{app:exp-details} specifies the training and evaluation protocol. Appendices~\ref{app:additional}--\ref{app:width} give settings for the additional experiments.


\item {\bf Experiment statistical significance}
    \item[] Question: Does the paper report error bars suitably and correctly defined or other appropriate information about the statistical significance of the experiments?
    \item[] Answer: \answerYes{}
    \item[] Justification: Tables report ten-seed means and population standard deviations. Appendix~\ref{app:paired} specifies paired tests and the absence of multiplicity correction.


\item {\bf Experiments compute resources}
    \item[] Question: For each experiment, does the paper provide sufficient information on the computer resources (type of compute workers, memory, time of execution) needed to reproduce the experiments?
    \item[] Answer: \answerNo{}
    \item[] Justification: The hardware and per-run times are reported in Appendix~\ref{app:exp-details}. A complete project-wide compute total is not reported.


\item {\bf Code of ethics}
    \item[] Question: Does the research conducted in the paper conform, in every respect, with the NeurIPS Code of Ethics \url{https://neurips.cc/public/EthicsGuidelines}?
    \item[] Answer: \answerYes{}
    \item[] Justification: The research conforms with the NeurIPS Code of Ethics.


\item {\bf Broader impacts}
    \item[] Question: Does the paper discuss both potential positive societal impacts and negative societal impacts of the work performed?
    \item[] Answer: \answerYes{}
    \item[] Justification: The introduction describes potential applications. Section~\ref{sec:conclusion} discusses the human and financial costs of inaccurate tail estimates in deployment.


\item {\bf Safeguards}
    \item[] Question: Does the paper describe safeguards that have been put in place for responsible release of data or models that have a high risk for misuse (e.g., pre-trained language models, image generators, or scraped datasets)?
    \item[] Answer: \answerNA{}
    \item[] Justification: The paper does not release data or models with a misuse risk.


\item {\bf Licenses for existing assets}
    \item[] Question: Are the creators or original owners of assets (e.g., code, data, models), used in the paper, properly credited and are the license and terms of use explicitly mentioned and properly respected?
    \item[] Answer: \answerNo{}
    \item[] Justification: The baseline papers are cited, but licenses for all implementation components are not documented in this manuscript.


\item {\bf New assets}
    \item[] Question: Are new assets introduced in the paper well documented and is the documentation provided alongside the assets?
    \item[] Answer: \answerNA{}
    \item[] Justification: The paper does not release new assets.


\item {\bf Crowdsourcing and research with human subjects}
    \item[] Question: For crowdsourcing experiments and research with human subjects, does the paper include the full text of instructions given to participants and screenshots, if applicable, as well as details about compensation (if any)? 
    \item[] Answer: \answerNA{}
    \item[] Justification: The paper does not involve human subjects.


\item {\bf Institutional review board (IRB) approvals or equivalent for research with human subjects}
    \item[] Question: Does the paper describe potential risks incurred by study participants, whether such risks were disclosed to the subjects, and whether Institutional Review Board (IRB) approvals (or an equivalent approval/review based on the requirements of your country or institution) were obtained?
    \item[] Answer: \answerNA{}
    \item[] Justification: The paper does not involve human subjects.


\pagebreak[3]
\item {\bf Declaration of LLM usage}
    \item[] Question: Does the paper describe the usage of LLMs if it is an important, original, or non-standard component of the core methods in this research? Note that if the LLM is used only for writing, editing, or formatting purposes and does \emph{not} impact the core methodology, scientific rigor, or originality of the research, declaration is not required.
    \item[] Answer: \answerNA{}
    \item[] Justification: The core method development in this research does not involve LLMs as any important, original, or non-standard components.

\end{enumerate}

\end{document}